# Multi-task Combinatorial Optimization: Adaptive Multi-modality Knowledge Transfer by an Explicit Inter-task Distance

**Peng Li** LIPENG@AMSS.AC.CN

**Bo Liu*** BLIU@AMSS.AC.CN

*Academy of Mathematics and Systems Science*
*Chinese Academy of Sciences, Beijing, China*

**Abstract**

Scheduling problems are often tackled independently, and rarely solved by leveraging the commonalities across problems. Lack of awareness of this inter-task similarity could impede the search efficacy. A quantifiable relationship between scheduling problems is to-date rather unclear, how to leverage it in combinatorial optimization remains largely unknown, and its effects on search are also undeterminable. This paper addresses these hard questions by delving into quantifiable useful inter-task relationships and, through leveraging the explicit relationship, presenting a speed-up algorithm. After deriving an analytical inter-task distance metric to quantitatively reveal latent similarity across scheduling problems, an adaptive transfer of multi-modality knowledge is devised to promptly adjust the transfer in forms of explicit and implicit knowledge in response to heterogeneity in the inter-task discrepancy. For faintly related problems with disappearing dependences, a problem transformation function is suggested with a matching-feature-based greedy policy, and the function projects faintly related problems into a latent space where these problems gain similarity in a way that creates search speed-ups. Finally, a multi-task scatter search combinatorial algorithm is formed and a large-scale multi-task benchmark is generated serving the purposes of validation. That the algorithm exhibits dramatic speed-ups of 2~3 orders of magnitude, as compared to direct problem solving in strongly related problems and 3 times faster in weakly related ones, suggests leveraging commonality across problems could be successful.

## 1. Introduction

Real-practice sequencing and scheduling problems seldom occur in isolation and often are related to each other. Imagine a cloud computing situation where multiple self-contained problems are emerging, necessitating they be simultaneously addressed. In reality, concurrent or closed-related problems are commonly tackled in isolation with zero ground knowledge, or only by utilizing its easily accessed domain knowledge. Such emerging multiple self-contained problems are rarely solved by leveraging commonality-based synergy benefits to foster search.

Optimization methods for scheduling have gone far without explicitly acknowledging inter-task relationships. Although remarkable progress has been made by continued innovations in speeding up optimization methods, both exact and approximate methods tackle the multiplicity of emerging problems in isolation. Task siloing means optimizing a new problem can become a Sisyphean dilemma. Ignorance of commonalities across problems frequently induces repetitive algorithm-tuning work that is dependent on expert knowledge, and induces futile searches accompanied by time-expensive evaluations. The unawareness of inter-task similarity could impede search efficacy.

Solving NP-hard scheduling problems with less human intervention but with accelerated convergences, even exponential speed-ups, seems a worthy place to start. We, consequently, charted an al-ternative route by appreciating the quantifiable, useful inter-task relationships, perceiving in this inter-task dependency the potential to exploit common knowledge shared among multiple self-contained problems. An accelerated search could be expected by transferring these commonalities from the source problem to the target problem, and thereby possibly alleviate time-consuming repeated searches. The implementation of knowledge transfer to accelerate search is largely based on the premise of quantitatively determining the distance between problems.

While the quantifiable relationship between scheduling problems is unclear, with the effects on search largely unknown, the relationship is non-trivial. Characterizing the relationship is an intricate task since the difference between combinatorial problems' landscapes can hardly be captured, mainly due to the curse of dimensionality. This paper attempts to answer computational challenges in discrepancy between problems, and to present an algorithm for leveraging an explicit relationship. Here, the algorithm is a collection of computationally found strategies specifying which problem can supply useful information to another, and by how much.

On the way to this destination, we take up four challenging questions. Can an explicit inter-task distance metric be derived to quantitatively reveal the latent similarity shared across scheduling problems? Can we guarantee to measure the distance in a sufficiently small overhead? Can an adaptive knowledge transfer scheme be designed under the guidance of the explicit inter-task distance metric to achieve speed-ups for correlated problems between which there is commonality? For weakly correlated problems between which there is only very faint similarity, are speed-ups in search achievable?

Our contributions are five-fold, and Figure 1 illustrates our approach.

For answering the first two questions, we constructed a normalized, symmetrical distance metric between combinatorial optimization problems (COPs) and focused on one of most intensively studied COPs, the permutation flowshop scheduling problem (PFSP). First, we found two new properties, i.e., the scale- and shift-invariance, and defined the order-isomorphism among problems. Second, we transformed the between-problem distance into the minimum distance from one problem to the order-isomorphic set of another problem, and followed by formulating a convex quadratic programming with an analytical solution. Then, by solving an edge-weighted bipartite matching problem, the distance metric was generalized to problems with differing dimensionality and with differently numbered jobs. That distance metric has extremely low computational complexity and good generalization to other permutation-based COPs. To our knowledge, this is the first attempt to measure the distance between COPs.

For the third question, an adaptive transfer of multi-modality knowledge was devised to foster search. To enrich the diversity of transferred knowledge for improved efficacy, explicit knowledge (partial solutions and complete solutions) and a novel implicit knowledge (solution evolution) were exploited and exchanged. Because the absence of prior knowledge about inter-task similarity leads to blindness in knowledge transfer, an adaptive scheme aware of inter-task distance was designed to promptly adjust the transfer in response to the heterogeneity in inter-task discrepancy. In combinatorial optimization it is the first step, to our knowledge, toward using explicit distance to guide transfer.

For the last question, how to transform the less related problems into closer problems to achieve the speed-up, we know the existence of cross-problem dependences is the premise of knowledge transfer, but for faintly related problems with disappearing dependences, a forcible knowledge transfer may lead to search stagnation. We learned the problem transformation function by a matching-feature-based greedy policy, and projected the problems into a joint latent space where they gain similarity. Commonality is leveraged to enhance the search by virtue of the implicit and explicit transfer across problems. Solutions to the transformed problems are inversely mapped back to their original domains. The mathematic analysis guarantees the equivalence and invertible property of the transformation. The (inverse-) projections of problems are relatively novel and can be used as a vital method in multi-task combinatorial optimization for domain adaption.

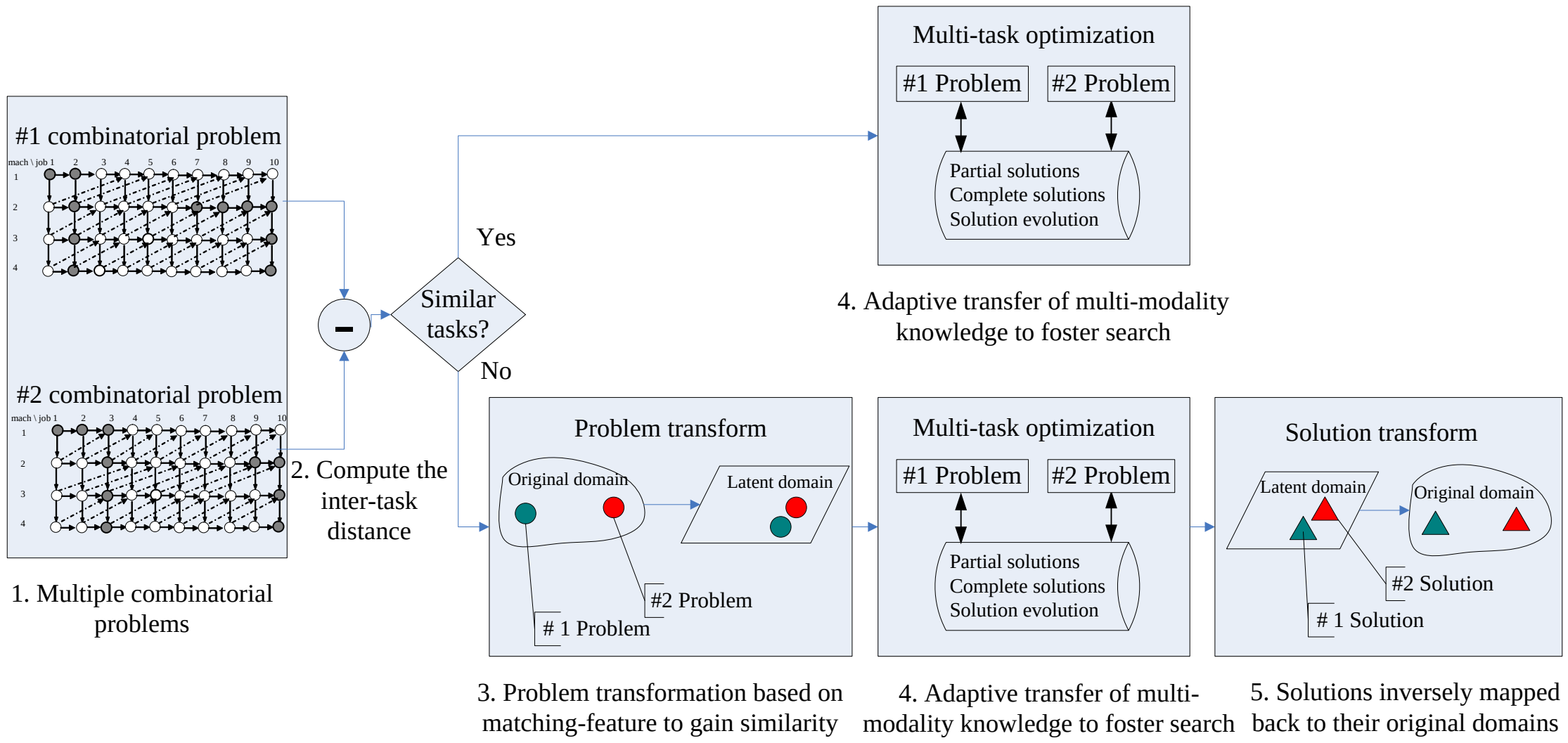

Figure 1: A schematic illustration of our approach. (1) Multiple combinatorial problems need to be solved. (2) Their similarity is explicitly computed by inter-task distance measure. (3) Problems with less similarity are projected into a latent space to increase commonality by the matching-feature-based problem transformation. (4) Adaptive multi-modality knowledge transfer is implemented to foster search, and inter-task distance is meant to control the transfer. (5) The solutions to the transformed problems are inversely mapped back to their original domains. Only Implement (4) for problems with high similarity.

We presented a multi-task combinatorial optimization algorithm that leverages the explicit relationship to accelerate convergences for solving multiple combinatorial problems. The algorithm explicitly computed the between-problem similarity by inter-task distance measure, increased the commonality for problems with less similarity by the matching-feature-based problem transformation, adaptively transferred multi-modality knowledge guided by inter-task distance, and performed solution perturbation globally and locally by scatter search.

For discerning the performance across different algorithms in the context of multitasking, a large-scale test bed is highly needed. However, there is no test bed for multi-task combinatorial optimization in the literature. Moreover, it was found that only faintly related problems were generated by random pairing instances from the well-known single-task Taillard's benchmark (Taillard, 1993). In this regard, we systematically generated a multi-task benchmark with 12,000 instances whose distances range from a value of one (completely unrelated) to zero (highly related). This is the first large-scale flowshop test bed in a multitasking context.

We verified through comprehensive experiments the efficacy of the multi-task combinatorial optimization algorithm. Remarkably, this approach exhibited dramatic speed-ups in strongly related problems but also in weakly related problems, suggesting its success lies in leveraging the commonality.

This research is vital to the field of machine learning for combinatorial problems where the distribution of learned samples remains a puzzle and generalization performances can only be known in hindsight. The inter-task distance metric, together with its usage in search, shed light on how to learn the underlying distribution of problem instances, how to determine whether the new case is an outlier or only a sample of the learned distribution, and how to leverage it.

## 2. Related Work

In scheduling, there is a belief that some common structure exists among different problems, a commonality that facilitates the reuse of solutions for solving new problems. However, since that important commonality across problems is not captured explicitly, it remains unknown. How to leverage it in combinatorial optimization remains largely unknown, and its effects on search are also undeterminable. Our research attempts to find this interaction and investigate its effects on search. Our research is related to a breadth of topics including: optimization methods for scheduling, case-based reasoning, problem reduction, multi-task Bayesian optimization, evolutionary multitasking, distance measures, and machine learning. Our concise review of topics most pertinent to scheduling algorithms which leverage the knowledge transfer, follows.

**Optimization methods for scheduling** evolved in the path of exact methods and approximate methods (dispatching rules, constructive heuristics, and improvement heuristics) (Beck, 2007; Liu, Wang, Liu, & Wang, 2011; Pinedo, 2012). Exact methods can guarantee the optimality but work best on problems with low- or middle- dimensionality. Approximate methods only provide feasible solutions on larger problems but cannot guarantee optimality, and improvement heuristics induce costly iterative searches. These methods tackle problems in isolation. No scheduling algorithms calculate explicit distance between problems to reveal the underlying commonality, let alone transfer commonality by a measurable distance. By contrast, we investigated the similarity measure between scheduling problems, and speed up the available method (we chose Scatter Search) by explicitly acknowledging these inter-task relationships.

**Case-based reasoning** retrieves the case(s) closer to the target case and builds a solution for the target case by reusing their solutions (Schmidt, 1998). The similarity between different scheduling problems is reflected in either the size of the intersection of jobs between source and target cases (Cunningham & Smyth, 1997; Kraay & Harker, 1997) or the Euclidean distance between fitness values of the solutions (Chang, Hsieh, & Wang, 2005). Nonetheless, surrogating distance measures do not provide a quantitative inter-task distance. This indicator, the overlap of jobs, cannot measure the whole picture alone. Due to the curse of dimensionality, it is intractable to characterize the inter-task differences by the differences in solutions' fitness values. To answer this based on mathematical analysis, we proposed an analytical inter-task distance metric by solving a constrained quadratic programming problem.

**Problem reduction** was originally meant to prove the computational complexity of a problem by transforming a generic instance of a problem to a particular instance of another. To obtain a concise theoretical analysis, the reductions used are relatively coarse-grained. Until recently, the problem re-duction was resorted to as a solver. Mathieson and Moscato (2020) in their study, proposed several reduction schemes mapping instances of Hamiltonian cycle problem to instances of traveling salesman problem (TSP). When afterward an effective TSP solver, the CONCORD (Applegate, Bixby, Chvatal, & Cook, 2006), was applied to solve TSP instances, the approach resulted in dramatic speedups. By contrast, we learned a problem transformation function by a matching-feature-based greedy policy, and then project problems into a joint latent space where they gain similarity.

**Multi-task Bayesian optimization** leverages the commonality between tasks to facilitate time-saving learning (Swersky, Snoek, & Adams, 2013). The expensive-to-evaluate function is represented by a cheap multi-task Gaussian surrogate model updated by feeding new observations. To lighten the burden of acquiring new observations, the inter-task dependence guides selection of promising points. This dependence is learned by a positive semi-definite covariance matrix over tasks (Swersky, et al., 2013), or by a transformation function (Bardenet, Brendel, Kégl, & Sebag, 2013). Though it sped up the search for continuous optimization (Pearce & Branke, 2018), it can not straightforwardly be extended to combinatorial optimization in that the covariance matrix of the Gaussian process becomes indefinite under combinatorial representations. Interestingly, a recent combinatorial Bayesian optimization using graph Cartesian product eliminates this defect (Oh, Tomczak, Gavves, & Welling, 2019), suggesting a crucial step towards the stage where multi-task Bayesian optimization can solve combinatorial problems.

**Evolutionary multitasking** inherits search behaviors from evolutionary algorithms and can solve multiple tasks at once by transferring knowledge via implicit chromosomes crossover (Bali, Ong, Gupta, & Tan, 2020; Gupta, Ong, & Feng, 2016) or in the explicit forms of solutions (Ding, Yang, Jin, & Chai, 2019; Feng et al., 2019). The few articles that thus far address vehicle routing (Feng et al., 2020), jobshop scheduling (Zhang, Mei, Nguyen, & Zhang, 2020), and TSP (Yuan, Ong, Gupta, Tan, & Xu, 2016), all assume the known inter-task relationships rather than explicitly computing them. For a prompt response to the heterogeneity in the inter-task discrepancy, we calculated explicit inter-task similarity and designed an adaptive transfer scheme that is aware of inter-task distance. Additionally, we mined the common knowledge shared by problems in novel forms, such as the partial solution (invariance partially revealing the common structure) and the solution evolution (iterative process of evolving to the best-so-far solution).

**Distance measures** have been largely investigated in fields of continuous optimization and machine learning, but quantifiable similarity between scheduling problems remains unknown. Because the difference between combinatorial landscapes can hardly be captured due to curse of dimensionality and none of smoothness (Barbulescu, Howe, Whitley, & Roberts, 2006; Streeter & Smith, 2006; Watson, Whitley, & Howe, 2005), characterizing the commonality is intricate. Available sampling-based methods like Spearman's correlation and maximum mean discrepancy (Gretton, Borgwardt, Rasch, Scholkopf, & Smola, 2012) are impracticable as their computational costs from sampling over the huge combinatorial search domain are enormous. Integral-based methods such as *KL* divergence (Bigi, 2003) fail to perform integrals over the discrete space without smoothness. To answer this deficiency, we regard the distance between two problems as the minimum distance from one problem to the order-isomorphic set of another problem. It can be formulated as a convex quadratic programming with an analytical solution.

**Machine learning**, leveraged to solve combinatorial problems (Aytug, Bhattacharyya, Koehler, & Snowdon, 1994; Bengio, Lodi, & Prouvost, 2021; Gratch & Chien, 1996; Hopfield & Tank, 1985;

Kool, Hoof, & Welling, 2019), learns the policies by imitation (supervised learning) or through experiences (reinforcement learning). Among the challenges machine learning encounters is that the true distribution of learned instances cannot be mathematically characterized in an explicit way. We do not, in fact, know how far the new instance is from the implicit distribution of the training set, suggesting generalization to unseen instances is a hard question. Deep learning excels when trained with a large data set but that situation is extremely scarce in combinatorial optimization. Recently, an effective graph convolutional neural network learned generalised policies for other stochastic shortest path problems by weight sharing between networks (Toyer, Thiebaux, Trevizan, & Xie, 2020). The inter-task distance metric we proposed could be a feasible tool for measuring distributions in machine learning to solve combinatorial problems.

The non-trivial relationship between scheduling problems is rather unclear and its effects on search largely unknown, albeit there is a belief that certain commonality exists among problems that facilitates search. This motivates us to answer the computational challenge in measuring similarity across scheduling problems, and to present an algorithm by taking advantage of the knowledge in common to foster search.

## 3. Permutation Flowshop Scheduling Problem

The intensively studied permutation flowshop scheduling problem (PFSP) is a notoriously intractable combinatorial optimization problem. It finds the permutation, say a sequence of jobs to be processed on machines, with respect to certain objective(s). Flowshop scheduling is an essential aspect of complex manufacturing and service facilities (Hall, 1998; Larsen & Pranzo, 2019; Liu, Wang, & Jin, 2007; Liu, et al., 2011).

In PFSP, a set of $n$ jobs $\{1, 2, \cdots, n\}$ has to be processed on each of the $m$ machines $\{1, 2, \cdots, m\}$. Each machine can execute at most one job at a time, and each job can be executed on at most one machine. The job permutation is kept the same on each machine. The $n$ jobs' permutation is denoted as $\pi = [\pi(1), \pi(2), \cdots, \pi(n)]$, where the $i$-th element $\pi(i)$, $i \in \{1,...,n\}$ is the job in the $i$-th position of the permutation. The processing time of job $\pi(i)$ on machine $j$ is given as $p_{\pi(i),j}$. The completion time for job $\pi(i)$ on machine $j$ is denoted as $C(\pi(i), j)$. The objective is to find a permutation $\pi$ to minimize the maximum completion time for all jobs on all machines, i.e., $C(\pi(n), m)$. The maximum completion time (makespan) is computed via the recursive Eqs. (2) - (4) (Liu, et al., 2007). Indices, parameters and variables used in this paper are shown in Appendix A. Nomenclature.

$$C(\pi(1),1) = p_{\pi(1),1} \tag{1}$$

$$C(\pi(i),1) = p_{\pi(i),1} + C(\pi(i-1),1), i = 2, \cdots, n \tag{2}$$

$$C(\pi(1), j) = p_{\pi(1),j} + C(\pi(1), j-1), \ j = 2, \cdots, m \tag{3}$$

$$C(\pi(i), j) = p_{\pi(i),j} + \max[C(\pi(i-1), j), C(\pi(i), j-1)], i = 2, \cdots, n \text{ and } j = 2, \cdots, m \tag{4}$$

$$f(\pi) = C(\pi(n), m) \tag{5}$$

To consistently represent problems with different scales and attributes, we give PFSP as a triplet $\{f, P, X\}$, where $f$ is the makespan, $P$ is the $(n \times m)$ problem specification matrix with $p_{i,j}$ as its element, and $X$ is the $(n \times n)$ solution matrix whose element is defined as

$$x_{i,r} = \begin{cases} 1, & if\ \pi(i) = r \\ 0, & otherwise. \end{cases} \tag{6}$$

We define $\varphi : \pi \to X$ as the mapping function from permutation $\pi$ to solution matrix $X$ in Eq. (6), and $\varphi^{-1} : X \to \pi$ as the inverse mapping function from $X$ to $\pi$ in Eq. (7).

$$\pi = [1, 2, \cdots, n] \cdot X^T \tag{7}$$

where $X^T$ is the transpose of solution matrix $X$, and $\pi$ and $X$ are equivalent.

## 4. New Properties for PFSP and Their Generalization

We establish for PFSP the scale- and shift- invariance in theorems 1 and 2, the theoretical basis for the inter-task distance measure in Section 5. We find an equivalent invertible transformation between PFSPs in theorems 3 and 4, these theorems being the support for the problem transformation to gain similarity in Section 7.

### 4.1 Scale- and Shift- Invariance for PFSP with Makespan

**Definition 1**: Problems $f$ and $g$ are order-isomorphic when $f(\pi_1) \le f(\pi_2) \Leftrightarrow g(\pi_1) \le g(\pi_2)$ for any solutions $\pi_1$ and $\pi_2$. The order-isomorphism between problems is represented as $f \cong g$.

**Theorem 1**: For any two problem specification matrices $P$ and $P'$, and for any positive scale value $t > 0$, If $P' = t \cdot P$, then for an arbitrary solution $\pi$, $f_{P'}(\pi) = t \cdot f_P(\pi)$ is true, where $f_P(\pi)$ and $f_{P'}(\pi)$ denote the makespans with $P$ and $P'$ under solution $\pi$, respectively.

Theorem 1 states that the PFSP possesses the scale-invariance property in terms of the makespan. By Definition 1, $f_{P'}$ and $f_P$ are order-isomorphic. The mathematical proof is in Appendix B. Proof of Theorem 1.

**Theorem 2**: For any two problem specification matrices $P$ and $P'$, If $P' = P + b \cdot E$, where $b \cdot E$ is a shift matrix $b \in R$, $R$ is the real number set and $E$ is an $(n \times m)$ matrix with all elements equal to 1, then $f_{P'}(\pi) = f_P(\pi) + (m + n - 1) \cdot b$ holds for any arbitrary solution $\pi$.

Theorem 2 states that PFSP possesses the shift-invariance property. By Definition 1, $f_{P'}$ and $f_P$ are order-isomorphic. Its proof is in Appendix C. Proof of Theorem 2.

Based on the theorems 1 and 2, the set for order-isomorphic problems can be defined.

**Definition 2**: An order-isomorphic problems set for function $f_P$ is defined as $G_p = \{ f_{P'} \mid P' = t \cdot P + b \cdot E,\ \ t > 0, b \in R \}$.

By Definition 2, the problems are kept order-isomorphic after the scaled and/or the shifted operations on $P$. The order-isomorphic problem set $G_p$ contains all problems generated by performing scaled or/and shifted operations on $P$. Geometrically, $P$ could be represented as a point in the $R^{n \times m}$ space, while $G_p$ is a family of rays passing $P$ with different intercepts in the $R^{n \times m}$ space. The distance from $f_Q$ to $f_P$ could be the minimum distance from $Q$ to $G_p$. We will elaborate on it later in this paper.

### 4.2 Equivalent Invertible Transformation between PFSPs

**Theorem 3**: For any two problem specification matrices $P$ and $P'$, if $P'=O\cdot P$, then $f_{P'}(\varphi^{-1}(X))=f_P(\varphi^{-1}(X\cdot O))$ holds for any solution matrix $X$.

The $\varphi^{-1}(\cdot)$ is defined in Eq. (7), denoting the mapping from solution matrix to a permutation. Matrix $O$ is an arbitrary $(n\times n)$ permutation matrix where each row and column of $O$ contain a single element 1, and the remaining elements are zero. The permutation matrix $O$ is referred to as the *transformation function*. $P$ is transformed to $P'$ by pre-multiplying $P$ by the transformation function $O$, that is, by swapping rows of $P$. Meanwhile, the solutions to $f_{P'}$ will change accordingly by performing the same transformation. Its proof is in Appendix D. Proof of Theorem 3.

**Theorem 4**: For any two problem specification matrices $P$ and $P'$, if $P=O^{-1}\cdot P'$ where $O^{-1}$ is the inverse of $O$, then $f_P(\varphi^{-1}(X))=f_{P'}(\varphi^{-1}(X\cdot O^{-1}))$ holds for any solution matrix $X$.

We refer to the permutation matrix $O^{-1}$ as the *inverse transformation function*. It transforms $P'$ back to $P$ by pre-multiplying $P'$ by the inverse transformation function. Meanwhile, the solutions to $f_P$ will change accordingly by performing the same transformation. Its proof is in Appendix E. Proof of Theorem 4.

By Theorems 3 and 4, the equivalent invertible transformation between PFSPs is achieved by the transformation function $O$ and its inverse transformation function $O^{-1}$. The bilateral transformations between $P$ and $P'$ are expressed by

$$P'=O\cdot P\Leftrightarrow P=O^{-1}\cdot P'. \tag{8}$$

### 4.3 Generalization to PFSP with Other Objectives

Although Theorems 1-4 are established for PFSP with makespan, they hold for PFSP with other objectives. We explain their generality for PFSP with continuous and discrete objectives, respectively.

The PFSP with respect to minimization of total completion time is a continuous function of processing time. Its objective is formulated

$$f(\pi)=\sum_{i=1}^{n}C(\pi(i),m). \tag{9}$$

Theorems 1, 3 and 4 are applicable. Theorem 2 changes to

$$f_{P'}(\pi)=f_P(\pi)+[n\cdot m+n\cdot(n-1)/2]\cdot b. \tag{10}$$

The PFSP with the goal of minimizing the number of tardy jobs is a discrete function of the processing time. A tardy job means its completion occurs after its due date. Its objective is formulated

$$f(\pi)=\sum_{i=1}^{n}U_i \tag{11}$$

$$U_i=\begin{cases}1 & if\ C(\pi(i),m)>d(\pi(i))\\ 0 & otherwise\end{cases} \tag{12}$$

where $d(\pi(i))$ is the due date for job $\pi(i)$. Theorems 1-4 hold by incorporating the due date into the problem specification matrix $P$ to formulate a new problem specification matrix as $[P,d]$. $E$ is an $n\times(m+1)$ matrix with all elements set to be one.

## 5. Inter-Task Distance Metric between PFSPs

In this section, challenges in quantitatively measuring the similarity between different PFSP problems are addressed. First, based on mathematical analysis, we propose a normalized, symmetrical inter-task distance metric by solving a constrained quadratic programming problem. Second, computational complexity analysis shows the distance can be computed in a sufficiently small overhead. Third, we generalize to problems with differing dimensionality by augmentation of problem dimensionality, and generalize to problems with differently numbered jobs by approximately solving an edge-weighted bipartite matching problem.

### 5.1 Constrained Quadratic Programming to Search for the Distance

We presume the distance $d(f_Q, f_P)$ from $f_Q$ to $f_P$ is the minimum distance from $f_Q$ to the order-isomorphic problems set $G_p$ of $f_P$, and regard the distance as residual error (i.e., the difference) between matrix $Q$ and its optimal approximation $P'$ in the space $G_p$, which is formulated as constrained quadratic programming of

$$\begin{aligned} d(f_Q, f_P) &= \min_{f_{P'}\in G_P} \|Q-P'\|_F \\ &= \min_{t>0,b}\|Q-t\cdot P-b\cdot E\|_F \end{aligned} \tag{13}$$

where $\|\,.\,\|_F$ is the Frobenius norm making Problem (13) convex, $f_{P'}$ is one function from the order-isomorphic problems set $G_p$, $P'$ is the problem specification matrix of $f_{P'}$, and $t$, $b$ and $E$ are defined in Theorems 1 and 2.

Problem (13) without the linear inequality constraint on $t$ is an unconstrained quadratic programming problem and can be optimally solved by the least square (Sorenson, 1970) with $t^0$ at its optimum. If $t^0$ is in the feasible region for Problem (13), it is the optimum for Problem (13), whereas if $t^0$ is not a feasible solution for Problem (13), it cannot be the optimum. Since $\|Q-t\cdot P-b\cdot E\|_F$ is continuous on $t$ and has lower bound, its infimum exists. The infimum could be the distance when $t=0$. Combining these two situations, the optimum $t^*$ can be explicitly represented as

$$t^* = \max[t^0, 0] \tag{14}$$

$$t^0 = \frac{n\cdot m\cdot\sum_{j=1}^{n}\sum_{i=1}^{m} q_{i,j}\cdot p_{i,j} - (\sum_{j=1}^{n}\sum_{i=1}^{m} q_{i,j})\cdot(\sum_{j=1}^{n}\sum_{i=1}^{m} p_{i,j})}{n\cdot m\cdot\sum_{j=1}^{n}\sum_{i=1}^{m} q_{i,j}{}^2 - (\sum_{j=1}^{n}\sum_{i=1}^{m} q_{i,j})^2}. \tag{15}$$

The optimum $b^*$ is

$$b^{*} = \frac{1}{n \cdot m}[\sum_{j=1}^{n}\sum_{i=1}^{m} q_{i,j} - t^{*} \cdot (\sum_{j=1}^{n}\sum_{i=1}^{m} p_{i,j})] \tag{16}$$

where $q_{i,j}$ and $p_{i,j}$ are elements of $Q$ and $P$, respectively.

By substituting $t^{*}$ and $b^{*}$ into (13), the distance is

$$d(f_Q, f_P) = \left\| Q - t^{*} \cdot P - b^{*} \cdot E \right\|_F . \tag{17}$$

Then, an analytical distance can be obtained as:

$$d(f_Q, f_P) = \left\| Q^{*} - t^{*} \cdot P^{*} \right\|_F \tag{18}$$

where

$$Q^{*} = Q - \frac{1}{n \cdot m} \cdot \sum_{j=1}^{n}\sum_{i=1}^{m} q_{i,j} \cdot E \tag{19}$$

$$P^{*} = P - \frac{1}{n \cdot m} \cdot \sum_{j=1}^{n}\sum_{i=1}^{m} p_{i,j} \cdot E . \tag{20}$$

### 5.2 Normalization of the Inter-task Distance

In line with the Cauchy-Buniakowsky-Schwarz inequality, the following inequality is satisfied

$$\left\| Q^{*} - t^{*} P^{*} \right\|_F^2 \geq \left( \left\| Q^{*} \right\|_F - t^{*} \left\| P^{*} \right\|_F \right)^2 \tag{21}$$

Then a normalized distance is

$$d(f_Q, f_P) = \begin{cases} \dfrac{\left\| Q^{*} \right\|_F - t^{*} \left\| P^{*} \right\|_F}{\left\| Q^{*} - t^{*} P^{*} \right\|_F}, & \left\| Q^{*} - t^{*} P^{*} \right\|_F \neq 0 \\ 0, & otherwise \end{cases} . \tag{22}$$

Distance value of 1 suggests completely unrelated PFSPs (i.e., PFSPs shares no similarities), while value of 0 suggests identical PFSPs, i.e., that there exists order-isomorphic relationship between PFSPs. Normalization of the inter-task distance is geometrically illustrated in Appendix F. Normalization of the Inter-task Distance.

### 5.3 Computational Complexity

The inter-task distance is expressed as a closed-form analytical solution to a convex quadratic programming problem. Its time complexity is extremely low, only in $O(m \cdot n)$. Detailed analysis is in Appendix G. Computational Complexity of Inter-task Distance Metric. It is observed that components involved in computing the inter-task distance are mainly matrix element operations, rather than matrix operations. Being free from the matrix structure, the inter-task distance can be calculated via parallel computing to further speed up, which ascertains a possible advantage in solving high-dimensional and real-time problems.

### 5.4 Generalize the Inter-task Distance through Preprocessing

The inter-task distance measure makes intuitive sense when the problems have the same dimensionality and jobs are appropriately numbered. In a situation where two PFSPs involve exactly the same job and same machines, but the jobs are differently numbered, their problem specification matrix would look different with rows swapped around. The inter-task distance measure would indicate low similarity despite the fact that the two PFSPs are identical. In a situation where two PFSPs differ in their number of machines and jobs, their problem specification matrixes would have different dimensionality and the inter-task distance metric could not work since the precondition of Problem (13) is that $Q$ and $P$ have identical dimensionality.

To slack the strong precondition, we proposed preprocessing procedures to generalize the distance measure to problems with differing dimensionality and with differently numbered jobs.

#### 5.4.1 AUGMENTING PROBLEM DIMENSIONALITY

One natural way to remedy the differing dimensionality issue relies on transforming one problem into a new one either by reducing or by augmenting its dimensionality so that both problems have the same dimensionality after transformation.

For continuous optimization, it is usual to marginalize out or condition on the extra variables to reduce dimensionality. While the same idea transferred to combinatorial optimization will bundle a set of jobs as a virtual job, that reduction fails to guarantee unchanged optima (Pinedo, 2012).

To protect the optima against changing after the problem is changed, we resorted to padding zeros in the problem specification matrix of the smaller problem. The problem does not change when either virtual jobs with zero processing time or virtual machines with zero processing time are added.

An example to illustrate is: Let $Q$ ($n_Q$ jobs and $m_Q$ machines) and $P$ ($n_P$ jobs and $m_P$ machines) be two problem specification matrixes. Without losing generality, we assume $n_Q < n_P$.

When $m_Q < m_P$, $P$ is unchanged and $Q$ is augmented as follows

$$Q \equiv \left( \begin{array}{c|c} Q & O_1 \\ \hline O_2 & O_3 \end{array} \right) \tag{23}$$

where $O_1$, $O_2$ and $O_3$ are zero matrixes sized as $n_Q \times (m_P - m_Q)$, $(n_P - n_Q) \times m_Q$ and $(n_P - n_Q) \times (n_P - n_Q)$, respectively.

When $m_Q > m_P$, $P$ and $Q$ are augmented as follows

$$\begin{aligned} P &\equiv \left( P \mid O_4 \right) \\ Q &\equiv \left( \begin{array}{c} Q \\ \hline O_2 \end{array} \right) \end{aligned} \tag{24}$$

where $O_4$ is the zero matrix with size of $n_p \times (m_Q - m_P)$.

When $m_Q = m_P$, $P$ is unchanged and $Q$ is augmented as in Eq. (24).

We noticed that a unified representation was proposed where search spaces of different problems were consolidated into one unified space by random key encoding. Interested readers may refer to Gupta, et al. (2016).

#### 5.4.2 Edge-Weighted Bipartite Matching

To remedy the issue on differently numbered jobs, we determined the matching relationship between jobs from different problems. Intuitively, matrix $Q$ is transformed by swapping its rows (re-ordering jobs in the matrix) so that the matrix $Q$ is closer to matrix $P$. We formulated an optimization problem that minimized the distance between transformed $Q$ and $P$ by searching for an optimal order of jobs in matrix $Q$

$$\min_{X_Q} d(f_{X_Q \cdot Q}, f_P) \tag{25}$$

where $X_Q$ denotes an $(n \times n)$ permutation matrix with exactly one entry of 1 in each row and each column and 0s elsewhere.

The search space for (a perfect matching) Problem (25) is a factorial of $n$. To reduce searching complexity, we transform it into an edge-weighted bipartite matching problem (Matoušek & Gärtner, 2007). The row of the matrix (indexed by job number) is regarded as the variable which is represented by the vertex. The edge weight is represented by the similarity between variables. We denote the edge weight matrix by $S$, and its element $s_{i,j}$ is the Pearson's correlation coefficient between the $i$-th row in $Q$ and the $j$-th row in $P$. The edge-weighted bipartite matching problem is formulated

$$\max_{X_Q} tr(X_Q \cdot S). \tag{26}$$

The optima $X_Q^*$ for Problem (26) can be solved by the Hungarian algorithm in polynomial complexity (Matoušek & Gärtner, 2007). Though Problem (26) is not equivalent to Problem (25), it can parsimoniously yield a satisfactory match for the perfect matching Problem (25).

## 6. Adaptive Transfer of Multi-Modality Knowledge

In a fictitious situation in aircraft formation in which two pilots, the lead pilot (target problem) and wingman (source problem), cannot look out for each other, they share a faint dependency. A forcible radio communication (knowledge exchange) would expose them to the enemy and countervail both pilots' endeavor on functioning (approaching the optima). In the context of optimization, the transfer (content and timing) should be promptly adjusted in response to the heterogeneity in discrepancy between problems.

For the purpose of adjustment, we designed an adaptive transfer of multi-modality knowledge to foster search. To improve efficacy, the diversity of transferred knowledge is enriched, with novel implicit knowledge (solution evolution) and explicit knowledge (partial solutions and complete solutions) being exploited and exchanged, rather than only transferring the complete or partial solutions in literature. To mitigate the blindness in knowledge transfer, an adaptive scheme aware of the inter-task distance is designed to promptly adjust transfer in terms of contents and timing. Hereafter, the source problem means the knowledge provider while the target problem is the receiver. The roles of providers or receivers can be interchanged.

### 6.1 Solution Evolution

A novel implicit knowledge, i.e., solution evolution, is proposed. Since generating a new solution by perturbation of the former solution is equivalent to pre-multiplying the old solution matrix by a permutation matrix, the product of a series of permutation matrices models the iterative process of evolving to the best-so-far solution. We refer to the product as the *solution evolution*. To enhance the knowledge sharing, the solution evolution is transferred from the source problem to the target problem. The new trial solutions to the target are generated by learning the source problem's solution evolution. Figure 2 details the implementation. Its added benefits will be shown in experiments.

**BEGIN**

**Input**: The initial solution to the source problem $\pi_s^{in}$, the initial solution to the target problem $\pi_t^{in}$ and the best-so-far solution to the source problem $\pi_s^{best}$.

1): Map permutations $\pi_s^{in}$ and $\pi_s^{best}$ to matrices $X_s^{in}$ and $X_s^{best}$, respectively, according to Eq. (6), and we get

$$X_s^{in} = \varphi(\pi_s^{in}) \text{ and } X_s^{best} = \varphi(\pi_s^{best}),$$

where $X_s^{in}$ and $X_s^{best}$ are the permutation matrices.

2): Quantify the evolving process from the initial solution $X_s^{in}$ to the best-so-far solution $X_s^{best}$ to obtain the *solution evolution* $O_\pi$ as follows.

Since $X_s^{in}$ and $X_s^{best}$ are two permutation matrices, there exists a permutation matrix $O_\pi$ such that

$$X_s^{in} \cdot O_\pi = X_s^{best}.$$

Since the inverse of a permutation matrix exists, we have

$$O_\pi = (X_s^{in})^{-1} \cdot X_s^{best}.$$

For any permutation matrix, it holds that

$$(X_s^{in})^{-1} = (X_s^{in})^T.$$

Then, the solution evolution is

$$O_\pi = (X_s^{in})^T \cdot X_s^{best}$$

3): Generate the new trial solution by learning the source's solution evolution $O_\pi$

$$\pi_t^{new} = \pi_t^{in} \cdot O_\pi.$$

**Output** $\pi_t^{new}$ to the target problem.

**END**

Figure 2: Knowledge Form: Solution Evolution

### 6.2 Partial Solution

It is established that the precedence between any two jobs usually affects the schedule's performance even when the two jobs are not adjacent. In other words, when two problems have similar structure, it happens more than not that a pair of jobs maintains the same precedence in high-quality solutions for both problems. We call it the *invariance* that any two jobs keep their precedence in both source and target problems. This invariance partially reveals the structure they have in common.

Based on this observation, we put forward the partial solution as a knowledge form. The jobs involved in the invariance are used to construct the partial solution where the precedence of jobs is preserved. Then the partial solution is iteratively constructed to be a complete solution by interpolating the unallocated jobs. Figure 3 details the implementation.

**BEGIN**

**Input**: Best-so-far solutions $\pi_s^{best}$ and $\pi_t^{best}$ of the source problem and target problem.

1): Calculate the invariance index for each job by

$$H(i) = \sum_{j=1, j\neq i}^{n} h(i,j)/(n-1), \quad i = 1, 2, \cdots, n$$

where $H(i)$ is the invariance index for job $i$, and $h(i,j)$ is the indicative function. The value of $h(i,j)$ is 1 when job $i$ and job $j$ keep the same precedence in $\pi_s^{best}$ and $\pi_t^{best}$; otherwise the value of $h(i,j)$ is 0. An illustrative example is in the Appendix H. An Illustrative Example for Computing the Invariance Index.

2): Sort jobs in the descending order of their invariance indexes. The first $(n-PT)$ jobs belong to the invariance set, while the rest of $PT$ jobs belong to the non-invariance set. Remark: $PT$ is set as $\min(\max(\sum_{i=1}^{n} I(H(i) > 0.5),\ 4),\ n/2)$, where $I$ is the indicative function and $I = 1$ when $H(i) > 0.5$; otherwise, $I = 0$.

3): A partial solution $\pi^{new}$ is generated by removing the jobs in non-invariance set from $\pi_t^{best}$.

4): Repeat until $\pi^{new}$ becomes a complete solution.

4.1): Pick the job with the lowest invariance index in the non-invariance set.

4.2): Insert the selected job to any possible position of $\pi^{new}$, to generate a series of candidate partial solutions.

Remark: In this scheme, the job with lower invariance index is preferred to other jobs in the non-invariance set, since it brings larger perturbation which is beneficial at the initial stages of constructing solutions.

4.3): Evaluate their performances for the target problem.

4.4): The candidate partial solution with the minimum makespan is set to be $\pi^{new}$.

4.5): The selected job is deleted from the non-invariance set.

**Output** $\pi^{new}$ to the target problem.

Figure 3: Knowledge Form: Partial Solution

### 6.3 Complete Solution

Transfer of the complete solution is straightforward. Our hypothesis is that if both problems are highly related, the optimal or sub-optimal solution for the source problem can maintain its quality in the target problem. Therefore, the solution would be ready for reuse in the target problem. We will investigate the transferability of optimal/best-so-far solutions between problems.

### 6.4 Adaptive Transfer Scheme

We report in advance here an interesting phenomenon that we revealed in the later experiment section. The best-so-far solution to the source problem can be ranked with a very high probability among the top 10 solutions to the target problem, when their inter-task distance is within 0.1. Although the transfer of the high-quality solution leads to dramatic speed-ups in this case, in another extreme case one problem's high-quality solution is deteriorated to a random solution for another problem when their inter-task distance is between 0.9 and 1.0. Under such circumstances, the knowledge transfer leads not to speed-ups but to hindering search.

Ideally, the knowledge transfer should be automatically adjusted by the inter-task commonalities, but combinatorial optimization literature never explicitly measures the distance between combinatorial problems nor ever leverages the distance to guide transfer. Rather, the concept of distances is frequently used in hindsight to empirically explain algorithm efficiency for particular problems.

We propose an adaptive transfer scheme guided by the explicit inter-task distance to determine whether/what to transfer at runtime. As an optimal decision problem, this minimizes the absolute value of the difference between the ratio of cumulative votes against transfers to all votes, and the inter-task distance. It follows

$$\underset{p\in\{1,2\}}{\operatorname{argmin}}\left|\frac{v(2)+(2-p)}{v(1)+v(2)+1}-d\right| \tag{27}$$

where the voting vector $v$ has two elements, $v(1)$, the cumulative votes in favor of transfers, and $v(2)$, the cumulative votes against transfers, with $d$ is the inter-task distance. Here, with these three forms of problem knowledge, we use the same formula for each form to control whether to trigger that transfer. We represent the voting vector for solution evolution, partial solution, and complete solution by $v_e$, $v_p$, and $v_c$, respectively, with these three voting vectors initialized to zeros.

When $p=1$, it votes against transfers and updates $v=v+[0,\ 1]$. When $p=2$, it votes in favor of tranfers and updates $v=v+[1,\ 0]$. According to Eq. (27), the inter-task distance $d$ controls whether to transfer. For instance, when $d$ is small (related problems), the possibility to implement transfers is high. When $d$ is large (faintly related problems), the possibility to implement transfers is low. When $d$ is 1 (completely unrelated problems), transfer is prohibited.

## 7. Problem Transformation based on Matching-Feature to Gain Similarity

It happens in reality that some problems share few commonalities. In the later experiment section, we show that almost all instances from the well-known Taillard's benchmark (Taillard, 1993) are unrelat-ed or faintly related. The forcible knowledge transfer between problems that lack similarity may lead to search stagnation.

This section focuses on how to transform weakly related problems into highly related ones. Based on Theorems 3 and 4, there obviously exists an equivalent invertible transformation between PFSPs (Eq. 8). The next important concern is how to find problem transformation functions that make problems more similar.

We propose a matching-feature-based problem transformation to project weakly-related problems into a latent space where they gain similarity. As detailed in Figure 4, the inter-task distance metric has the role of measuring inter-task closeness after problem transformation. In particular, each problem is characterized by independent features (Step 2 in Figure 4), and we calculate correlation values between features of different problems to obtain the matching-feature relationship (Step 3). The matching-feature is used to construct the transformation functions $O_P$ and $O_Q$ (Step 4). The original pair of problems is transformed and jointly projected into a latent space by the transformation functions (Step 5). Computing their inter-task distance makes it easy to check whether they gain similarity in this new space or not (Step 6). If similarity is gained, the original pair of problems is replaced by the new pair (Step 7).

We adopt a greedy rule to build the transformation functions by matching-feature (see Step. 4 in Figure 4), that is, the most correlated pair of jobs has the highest priority to be enrolled. While the merit of the greedy rule is in its computational efficiency, its disadvantage is that the correlation between the latter pairs of features becomes weak. However, this defect has less impact on the optimization. We want to remark that there exists optimal rule whereby the matching-feature problem can be formulated as a combinatorial optimization problem (see Eq. 28). The optimal transformation functions $O_P$ and $O_Q$ can be obtained by solving it:

$$\min_{O_P, O_Q} d(f_{O_P \cdot P}, f_{O_Q \cdot Q}) \tag{28}$$

where the distance $d(f_{O_P \cdot P}, f_{O_Q \cdot Q})$ can be computed according to Eq. (13). To save computation budgets, we choose the greedy rule over the optimal rule.

**BEGIN**

**Input**: The original source problem $f_P$, target problem $f_Q$, and their inter-task distance $d(f_P, f_Q)$.

1) Initialize the transformation functions $O_P$ and $O_Q$ for $f_P$ and $f_Q$ as $(n \times n)$ null matrixes, respectively.
2) Decompose the $(n \times m)$ problem specification matrixes $P$ and $Q$ into $n$ row vectors $\{P_{(1)}, P_{(2)}, P_{(3)}, \cdots, P_{(n)}\}$ and $\{Q_{(1)}, Q_{(2)}, \cdots, Q_{(n)}\}$. Each vector corresponds to one job and represents its processing time on the $m$ machines. Hereafter, each vector is treated as a feature.
3) Calculate the Pearson's correlation coefficients between features of $\{P_{(1)}, P_{(2)}, P_{(3)}, \cdots, P_{(n)}\}$ and $\{Q_{(1)}, Q_{(2)}, \cdots, Q_{(n)}\}$ to build an $(n \times n)$ matching-feature matrix $MF = (mf_{i,j})$, which is constructed where $mf_{i,j}$ is the Pearson's correlation coefficient between $P_{(i)}$ and $Q_{(j)}$.
4) Compute the *transformation functions* $O_P$ and $O_Q$

    **For** $i = 1 : n$.

    4.1) Find the largest element in $MF$ and its indices $(idx_P,\ idx_Q)$.

    4.2) Set the elements in $idx_P$-th row and the $idx_Q$-th column to be null.

    4.3) Update transformation functions $O_P$ and $O_Q$ by setting $O_P(i, idx_P) = 1$ and $O_Q(i, idx_Q) = 1$.

    **End For**
5) Transform matrix $P$ and $Q$ into $P^{new}$ and $Q^{new}$ by

$$P^{new} = O_P \cdot P \text{ and } Q^{new} = O_Q \cdot Q.$$

6) Compute the inter-task distance $d(f_{P^{new}}, f_{Q^{new}})$.
7) If $d(f_{P^{new}}, f_{Q^{new}}) < d(f_P, f_Q)$, replace $P$, $Q$ by $P^{new}$ and $Q^{new}$. Otherwise, reset the transformation functions $O_P$ and $O_Q$ as identity matrix.

**Output** $f_P$ and $f_Q$, and $O_P$ and $O_Q$.

**END**

Figure 4: Matching-Feature-based Problem Transformation

## 8. Multi-task Combinatorial Optimization Algorithm

In previous chapters, we proposed the inter-task distance metric, adaptive transfer of multi-modality knowledge, and problem transformation based on matching-feature. In this section, we incorporate these techniques into a scatter search algorithm to build a multi-task combinatorial optimization algorithm.

### 8.1 Scatter Search

Scatter search was originally developed for solving integer programming, i.e., creating composite decision rules and surrogate constraints (Glover, 1977, 1998; Hvattum, Duarte, Glover, & Marti, 2013; Martí, Laguna, & Glover, 2005). Later, it was successfully applied to a wide range of hard combinatorial optimization problems (Martí, 2006), including scheduling (Burke, Curtois, Qu, & Vanden Berghe, 2010; Naderi & Ruiz, 2014; Yang, Li, Wang, Liu, & Luo, 2017).

Glover (1998) identified a template for implementation of scatter search, as shown in Figure 5 with subroutines in italics. We will detail their implementation for our problem.

**BEGIN**
  **Initialize**:
    Use *Diversification Generation* to generate a population of trial solutions;
    Use the *Reference Set Build Method* to build a reference set consisting of best solutions.
  **While (**Stopping conditions are not satisfied**)**
    Generate new subsets with the *Subset Generation* by operating on the reference set;
    Apply *Solution Combination Method* to a selected subset of solutions to obtain one or more new trial solutions;
    Apply *Improvement Method* to enhance the trial solutions;
    Apply *Reference Set Update Method* to update the reference set.
  **End While**
  Output Reference Set
**END**

Figure 5: The basic procedure of scatter search

### 8.2 Diversification Generation

For the purpose of the Diversification Generation, which is to generate an initial population of $N$ trial solutions with certain diversity and quality, we developed a composite generation rule. After using NEH (Nawaz, Jr, & Ham, 1983), CDS (Campbell, Dudek, & Smith, 1970), NEHKK1 and NEHKK2 (Dong, Huang, & Chen, 2008; Kalczynski & Kamburowski, 2007, 2008) to each generate a solution, the RAER (Ribas, Companys, & Tort-Martorell, 2010) creates the remaining solutions. It is expected that multiple diversified high-quality solutions could be obtained through these heuristics.

### 8.3 Reference Set Build and Update

At the beginning, a reference set is built by selecting the best $b$ solutions from the initial population of trial solutions. Then, along the search, the reference set is updated by maintaining $b$ "best" solutions found so far and ranking the reference solutions in the ascending order of their objective values.

### 8.4 Subset Generation

The Subset Generation generates a new subset by operating on the reference set. We adopt the simplest form of 2-element subset generation method, which generates all pairs of reference solutions.

## 8.5 Solution Combination

The Solution Combination generates a new trial solution by splicing the reference solutions in the subset. We adopt the solution combination method in Glover (1994). Two reference solutions are denoted as $\pi_1$ and $\pi_2$, where $\pi_1$ has a smaller makespan than $\pi_2$. Both reference solutions vote for the first unassigned position in the trial solution as follows. The first job in the trial solution is the first job in $\pi_1$, and then the successors of the last selected job in $\pi_1$ and $\pi_2$ become possible candidates for the trial solution. If the successor previously has been enrolled in the trial solution, its successor is selected. If the last selected job is at the end of the reference solution, its successor is the first job among the remaining jobs in the reference solution. If the two candidates are identical, any candidate is appended to the end of the trial solution. If they are different, they have identical probability to be selected to be appended to the end of the trial solution.

## 8.6 Improvement Method

The Improvement Method performs fine search by exploiting the promising region around the trial solution to achieve an enhanced search. In combinatorial optimization, simulated annealing, despite its simplicity, achieves state-of-the-art performances (Ball, Branke, & Meisel, 2018; Kirkpatrick, Gelatt, & Vecchi, 1983; Michiels, Aarts, & Korst, 2007). We utilize simulated annealing to improve the trial solution.

The simulated annealing starts from an initial state, that is, from the trial solution, and a new state is sampled from its neighbors. The Insert-based operation (Liu, et al., 2007) used to perform sampling, having randomly chosen two distinct elements in the trial solution, then inserts the back one before the front one. Different from the solution acceptance mechanism based on the survival of the fittest in the greedy search, the new candidate is accepted with a probability $P_T = \min\{1, \exp(-\Delta E / T)\}$, where $\Delta E$ is the change in the objective function value, and $T$ is the parameter mimicking the temperature in the annealing process. Worse candidate solution could be accepted in the probability which is gradually reduced as regards the annealing schedule.

The initial temperature $T_0$, annealing schedule, and steps of metropolis sampling significantly impact the simulated annealing's efficacy (Ball, et al., 2018). We use $T_0 = \sum_{i=1}^{n}\sum_{j=1}^{m} p_{i,j} / (10 \cdot n \cdot m)$ (Stutzle, 1998). The exponential annealing schedule used follows $T_k = \lambda \cdot T_{k-1}$, where cooling coefficient $\lambda = 0.9$, $T_k$ and $T_{k-1}$ are temperatures at generations $k$ and $k-1$, respectively (Liu, et al., 2007), and the steps of metropolis sampling are set to be $n \cdot (n-1)$.

## 8.7 Overall Framework of MTCO

We have sophisticatedly designed all the components prerequisite for instantiation of the multi-task combinatorial optimization algorithm (MTCO). Figure 6 illustrates the framework with a case in which the transfer is from the source problem to the target one. In MTCO, the transfer is bidirectional. The source problem is the knowledge provider, the target one is the receiver, and the roles of providers or receivers are interchanged along the search.

In MTCO, the distance between problems is measured by *inter-task distance metric*. Based on the quantifiable distance, the pair of problems is classified as either weakly correlated or strongly correlated. The decision boundary for classification is set to be 0.5 by experiments (see Section 9.8). Only weakly correlated problems are transformed into closely correlated problems by the *matching-feature-based problem transformation*. For each problem, an initial population of trial solutions is generated by the *diversification generation* method. The "best" $b$ solutions from the initial population build the reference set. New trial solutions are generated from two sources. Whereas one trial solution is generated by the *solution combination* following the *subset generation* procedure, once the transfer from the source to the target is activated, other trial solutions to the target are generated by implementing *adaptive transfer of multi-modality knowledge*. The best one among the new trial solutions is improved by the simulated annealing based *improvement method*. The *reference set is updated* by comparing the improved trial solution with the reference solutions. Repeat the above procedures until the termination condition is met. For transformed problems, the best solution in the reference set is mapped back to the original problem space by the *inverse transformation function*.

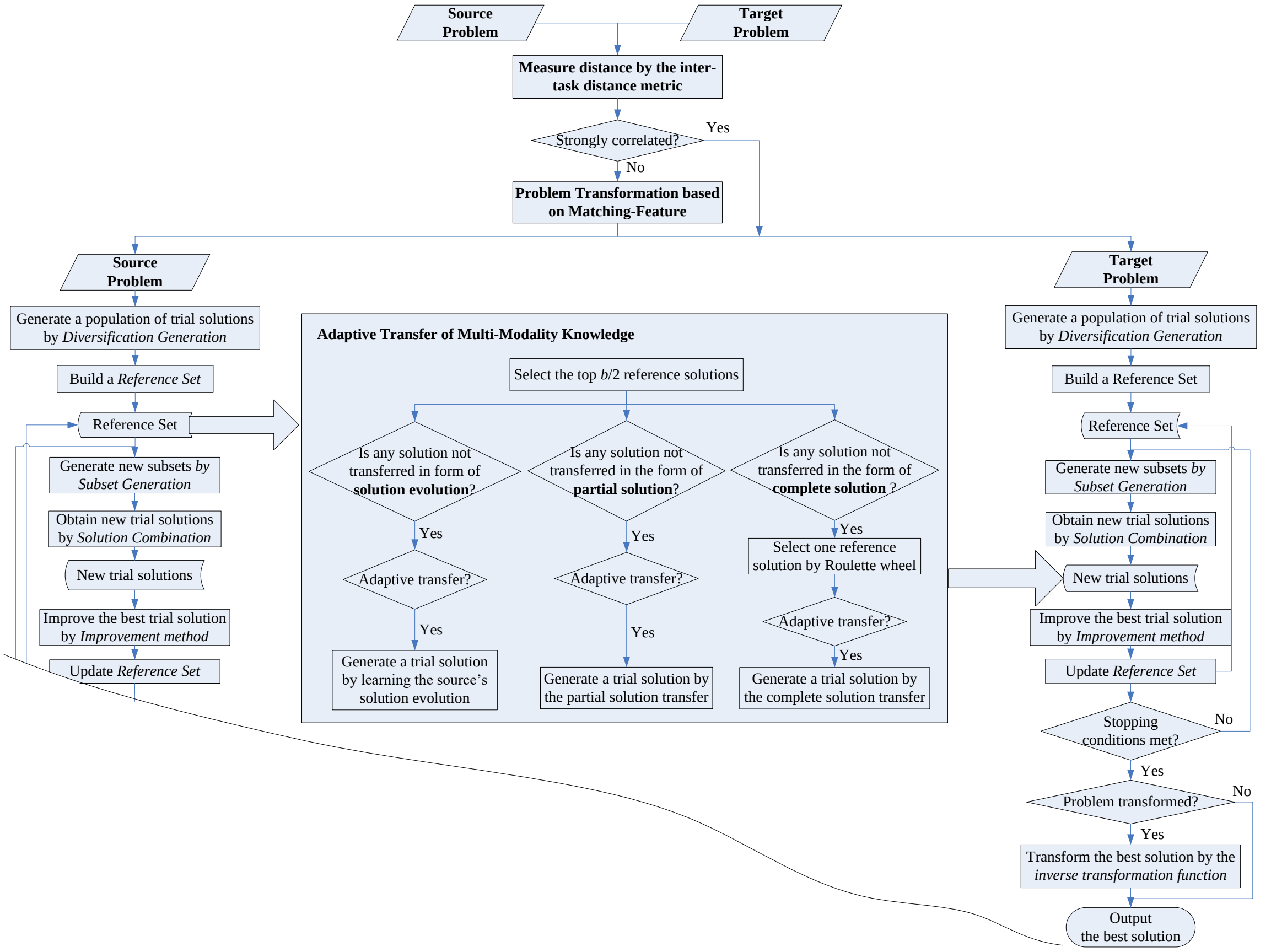

Figure 6: Overall framework of the multi-task combinatorial optimization algorithm

## 9. Benchmark, Test and Comparisons

In this section, we verify through comprehensive experiments the efficacy of the multi-task combinatorial optimization algorithm. For validation, we systematically generate a multi-task benchmark with 12,000 instances, and compare MTCO's performances against the single-task algorithm. Last, we analyze the impact on the searching performance of the adaptive transfer of multi-modality knowledge and the matching-feature based problem transformation.

## 9.1 Designing Multi-task Benchmark for PFSPs

### 9.1.1 INITIATIVES

Different from the single-task benchmark (i.e., when each instance contains only one optimization problem), each instance of the multi-task benchmark consists of at least two optimization problems. Since there is no test bed for multi-task combinatorial optimization, we must generate a benchmark to discern the performance difference across different algorithms.

Though it is straightforward to generate an instance by randomly choosing two instances from the available single-task benchmark, the following two experiments demonstrate that only instances of unrelated problems are generated. This simple way, called a random pairing method, failed to generate instances of highly correlated problems.

### 9.1.2 EXPERIMENT 1: RANDOM PAIRING METHOD

In the first experiment, we computed the distance between the optima of two instances from the well-known Taillard's benchmark (Taillard, 1993). We selected 90 instances whose optima are available in OR-library and, after pairing instances having equal dimensionality, we determined that 362 pairs (instances) can consolidate a tentative multi-task benchmark. The distance between optima for two problems within the same instance is calculated by the well-known precedence-based distance metric (Reeves & Yamada, 1998). Values of 0, 0.5, and 1 respectively indicate identical, unrelated, and totally inverse relationship between the optima of two problems. The detail computing procedure is in Appendix I. Precedence-based Distance Metric to Measure the Distance between Optima.

The frequency of observations (counts of instances) occurring in certain distance bins is depicted in Figure 7. The histogram reveals that for 334 instances out of the 362 (nearly 93%), their distances fall in the interval [0.4, 0.6]. For 231 of 362 instances (nearly 64%), their distances fall in the interval [0.45, 0.55]. Because most of the single-task instances are unrelated in terms of the precedence-based distance metric, this suggests the random pairing method fails to generate instances of highly correlated problems.

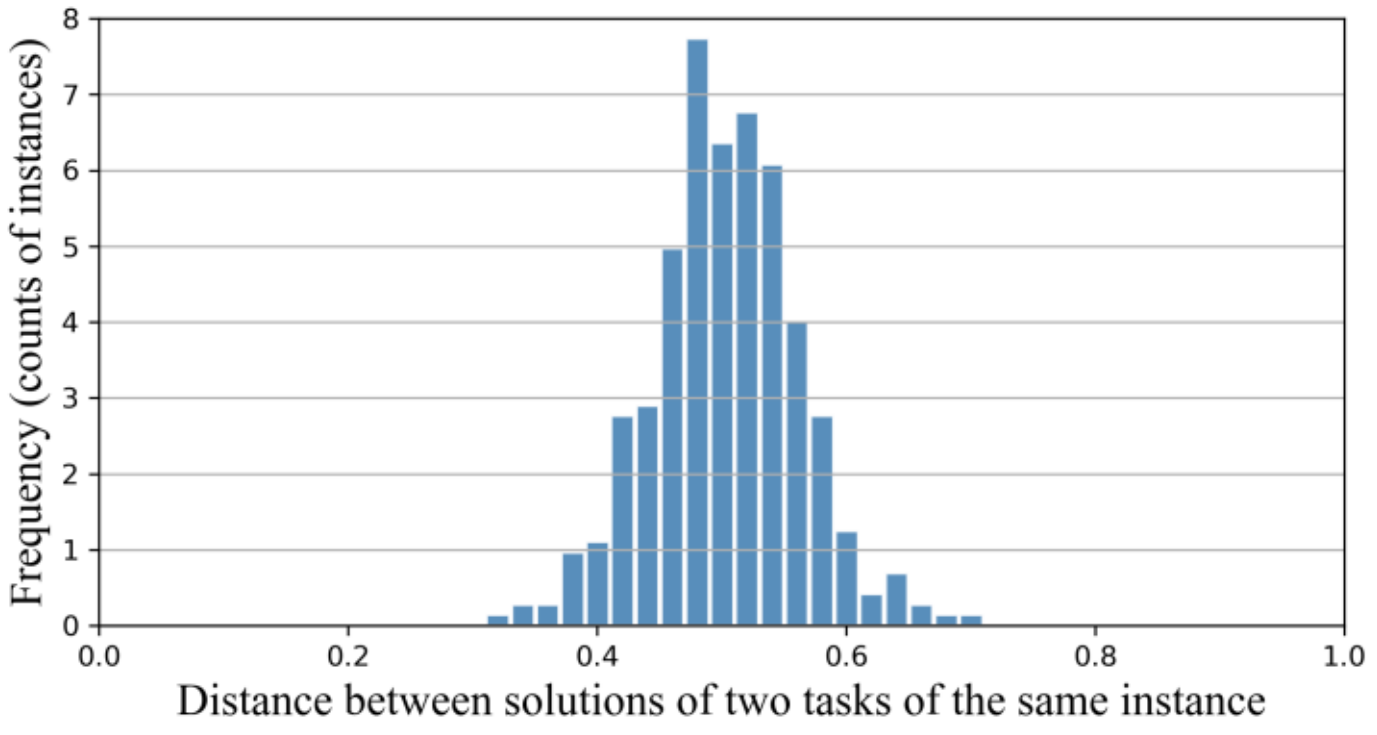

Figure 7: Frequency distribution histogram giving the frequency of observations (counts of instances) occurring in different distance bins.

#### 9.1.3 Experiment 2: Random Pairing Method

In the second experiment, we measured the distance between two problems using our proposed inter-task distance metric. By pairing Taillard's instances of identical dimensionality, multiple inter-task instances are generated and their inter-task distances are measured.

The results for all instances with [5, 10, 20] machines and [20, 50, 100, 200, 500] jobs are in Appendix J. Computational Results of Inter-task Distance for All Instances. The inter-task distances for instances with 50 jobs and 20 machines (i.e., pairwise combinations of Tai-51, Tai-52, …, Tai-60) are illustrated in Figure 8. The heatmap reveals that all of the inter-task distances fall in [0.94, 1], which suggests an unrelated relationship.

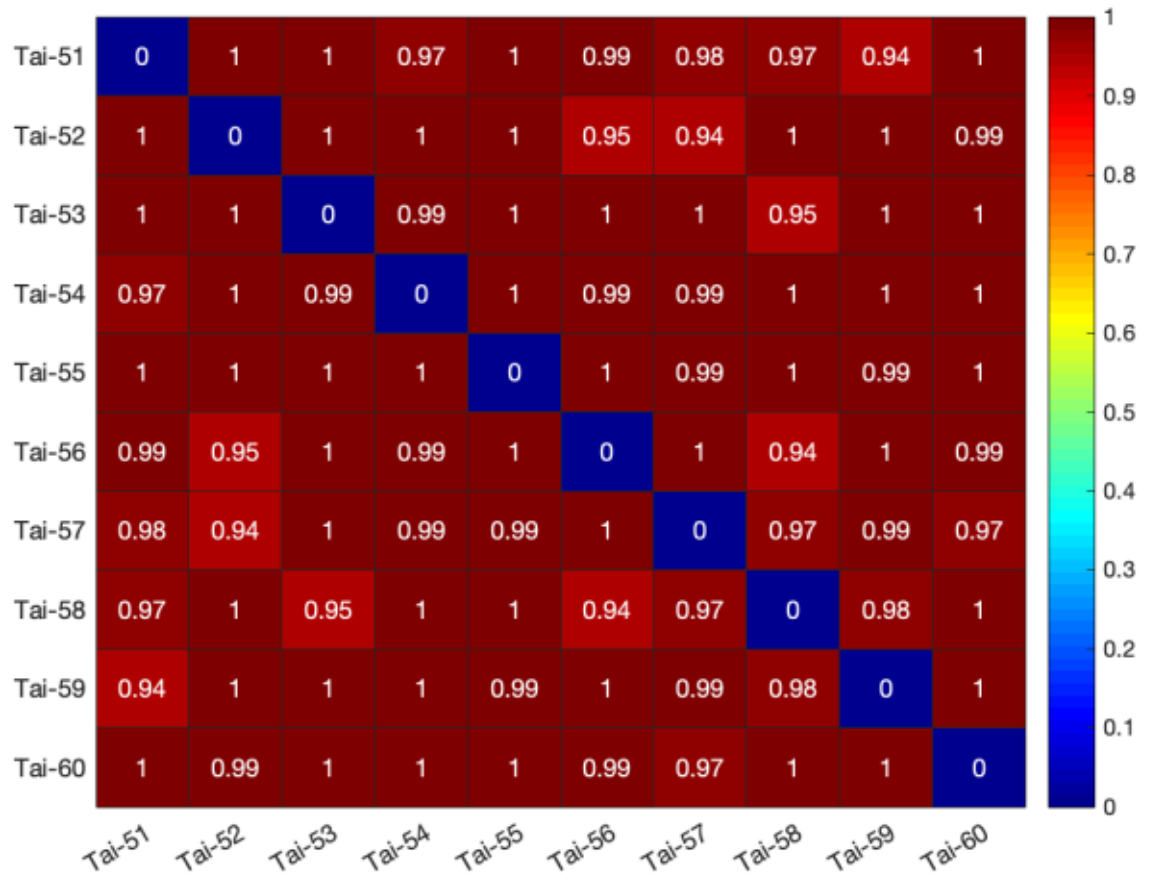

Figure 8: The inter-task distance instances' heatmap (50 jobs and 20 machines).

To sum up, almost all instances from the single-task benchmark are unrelated or faintly related in terms of the distance between optima as well as the inter-task distance metric. The phenomena motivate us to design the multi-task benchmark, rather than use the random pairing method.

#### 9.1.4 Probabilistic Element Replacement Method

To generate the multi-task benchmark, we proposed a probabilistic element replacement method. We illustrate the generation of a multi-task instance in Figure 9. Each multi-task instance consists of two problems, Problem 1 and Problem 2. Problem 1 is exactly the Taillard's instance. Problem 2 is generated via replacing elements of Problem 1's problem specification matrix by random numbers subjected to a uniform distribution $U(1,99)$ with a user-specified replacing probability $p_r$. To elaborate, for each element $p_{i,j}$ of Problem 2, if the random number $rand_{i,j}$ is less than $p_r$, then the element $p_{i,j}$ for Problem 2 is numbered by sampling $U(1,99)$; otherwise, $p_{i,j}$ adopts the value of corresponding element in Problem 1's problem specification matrix.

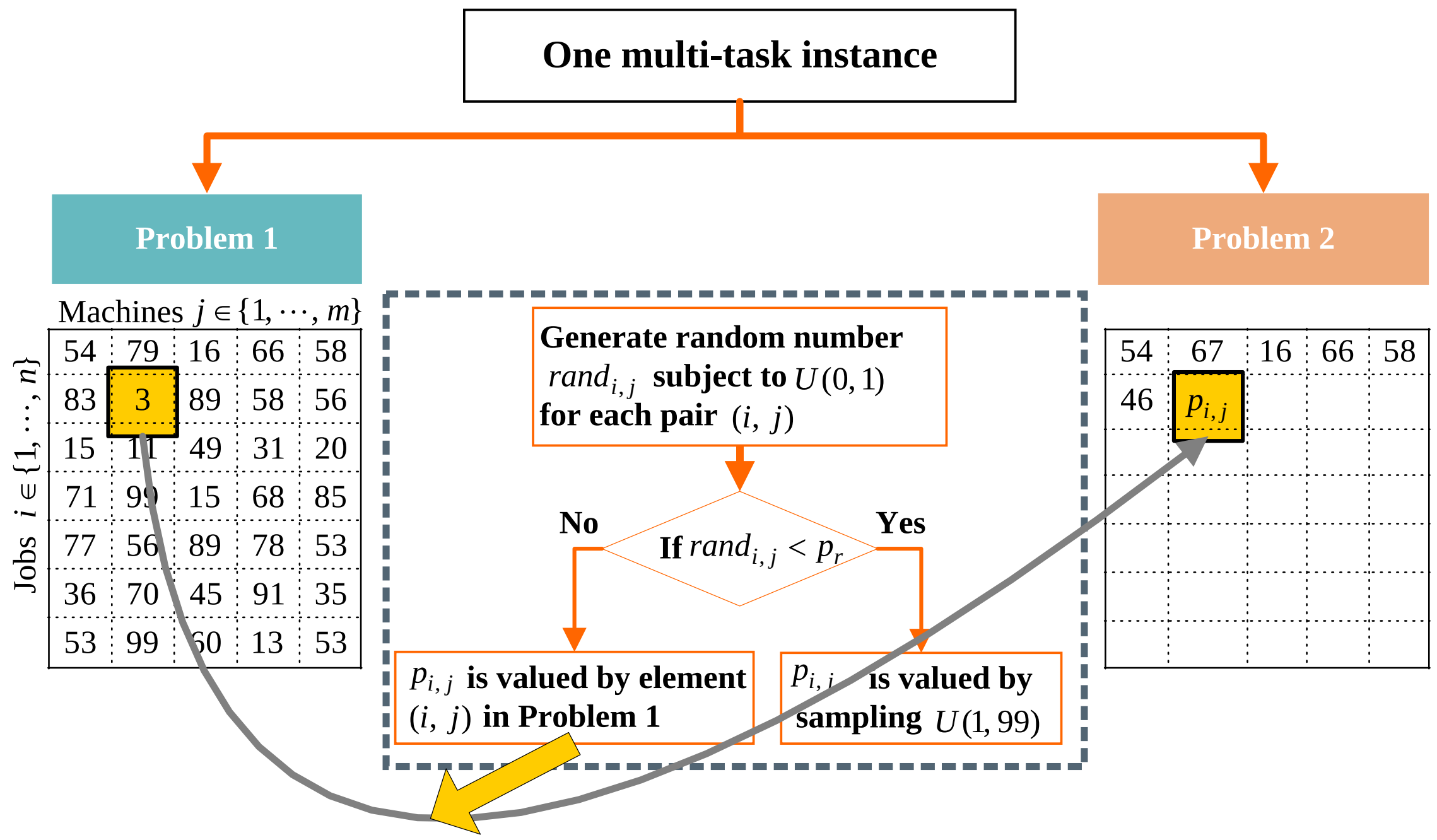

Figure 9: Illustration of generation of the multi-task instance.

### 9.1.5 EFFECTS OF REPLACING PROBABILITY ON DISTANCE DIVERSITY

The distance diversity is the core requirement in design of the multi-task benchmark. We investigate the impact of replacing probability $p_r$ (the only controlling parameter) on the distance variations.

In the Taillard's benchmark, there are total 120 instances, with $p_r$ ranging from (0, 1] with 100 equal intervals. In total, 12,000 instances are generated and their inter-task distances are measured. The inter-task distances are depicted against $p_r$ in Figure 10. The scatter plot reveals that the replacing probability is positively correlated with inter-task distance.

One special case is to set $p_r$ as 0, which means both problems are identical with distance of 0. Another special case is to set $p_r$ as 1, under which all $p_{i,j}$ for Problem 2 are sampled from $U(1,99)$. In this regard, Problem 2 could be regarded as an instance randomly selected from Taillard's benchmark, since the Taillard's instances are generated by sampling from $U(1,99)$. From Figure 10, it could be seen that the inter-task distance of instances generated with $p_r$ of 1 is about 1, which is consistent with the findings in Section 9.1.3, that instances from single-task benchmark are unrelated or faintly related in terms of the inter-task distance metric.

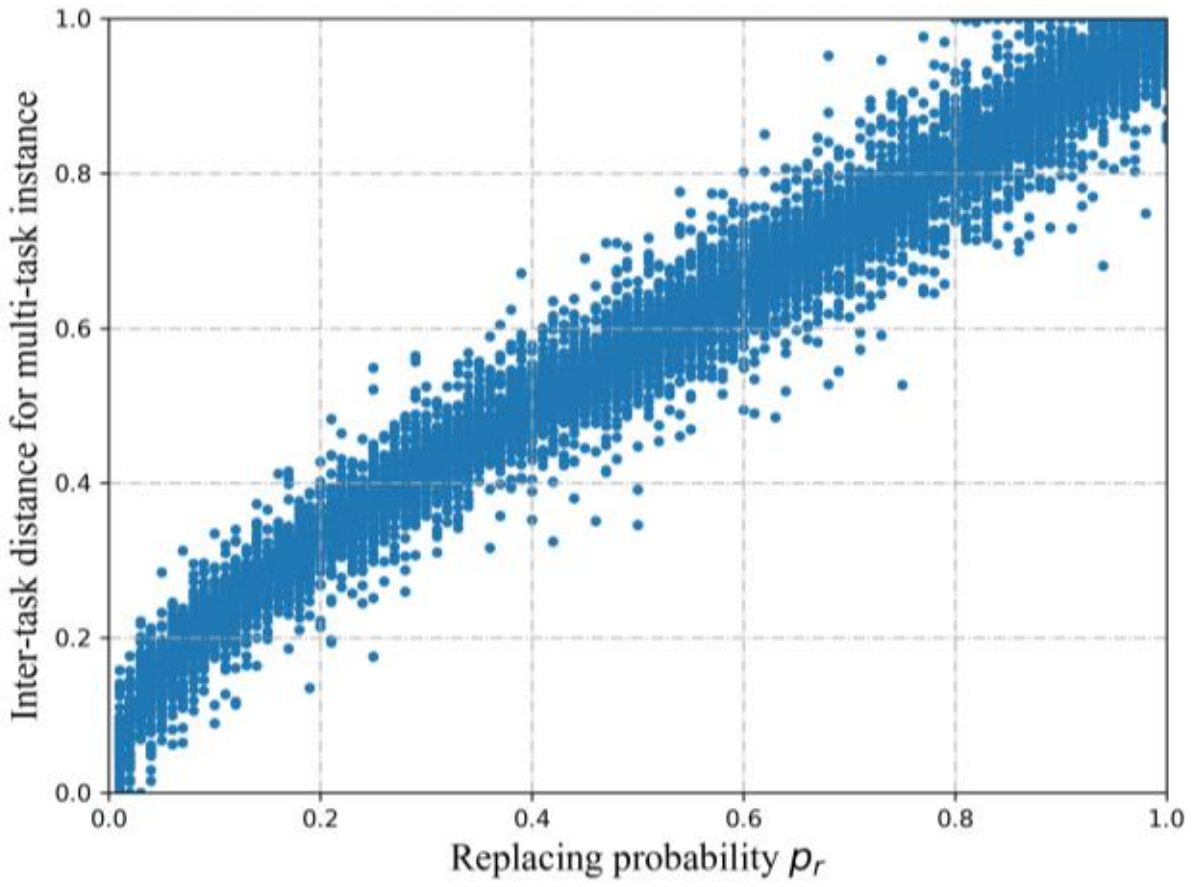

Figure 10: The effect of different replacing probabilities on the inter-task distance.

To sum up, we systematically generated a new multi-task benchmark. It contains a set of 12,000 instances and is publicly available in Li and Liu (2020). Through manipulating the single parameter $p_r$, the distance varies from one (completely unrelated problems) to zero (highly related problems).

## 9.2 Validating the Inter-Task Distance

We validate the correctness of the proposed inter-task distance measure. We investigate the consistency between the inter-task distance and the Spearman's rank correlation coefficient (SRCC).

As a nonparametric measure, SRCC assesses the monotonic relationship between two variables by computing the statistical dependence between the observations (i.e., rankings) of two variables (Mood, Graybill, & Boes, 1973). Intuitively, high correlation value suggests a similar rank between observations of two variables (a correlation of 1 for identical rank) while low value suggests a dissimilar rank (a correlation of −1 for fully opposed rank).

In our study, SRCC uses a moderate sampling size of 100,000 feasible solutions. The objective function Eq. (5) for both problems in a multi-task instance could be seen as two variables. For each variable, 100,000 observations could be obtained by evaluating the 100,000 feasible solutions through Eqs. (1)–(5). The raw numbers of the observations are converted into ranks within the variable. The statistical dependence between two variables is assessed by SRCC, the value of which is considered to be the similarity. We computed the SRCC value for each of 12,000 instances. The similarity between problems in each instance was also computed by inter-task distance measure.

To investigate the degree of consistency, the SRCC is depicted against inter-task distance measure with line of best fit in Figure 11. The scatter plot reveals that as the inter-task distance grows from 0 to 1, SRCC has an opposite trend, gradually decreasing from 1 to 0. Another obvious aspect seen in the scatter plot is the strong negative linear correlation with a goodness of fit of 0.9406 between these two measures, which suggests a very high degree of consistency between inter-task distance measure and SRCC. The similarity could be effectively computed, we conclude, by comparing the computational complexity of the inter-task distance measure to the SRCC measure. To be specific, SRCC's time complexity is in $O(Z \cdot m \cdot n + Z)$ where $Z$ is the sampling size, while inter-task distance meas-ure's time complexity is extremely low, only in $O(m \cdot n)$. The proposed inter-task distance metric is a faithful and efficient similarity metric.

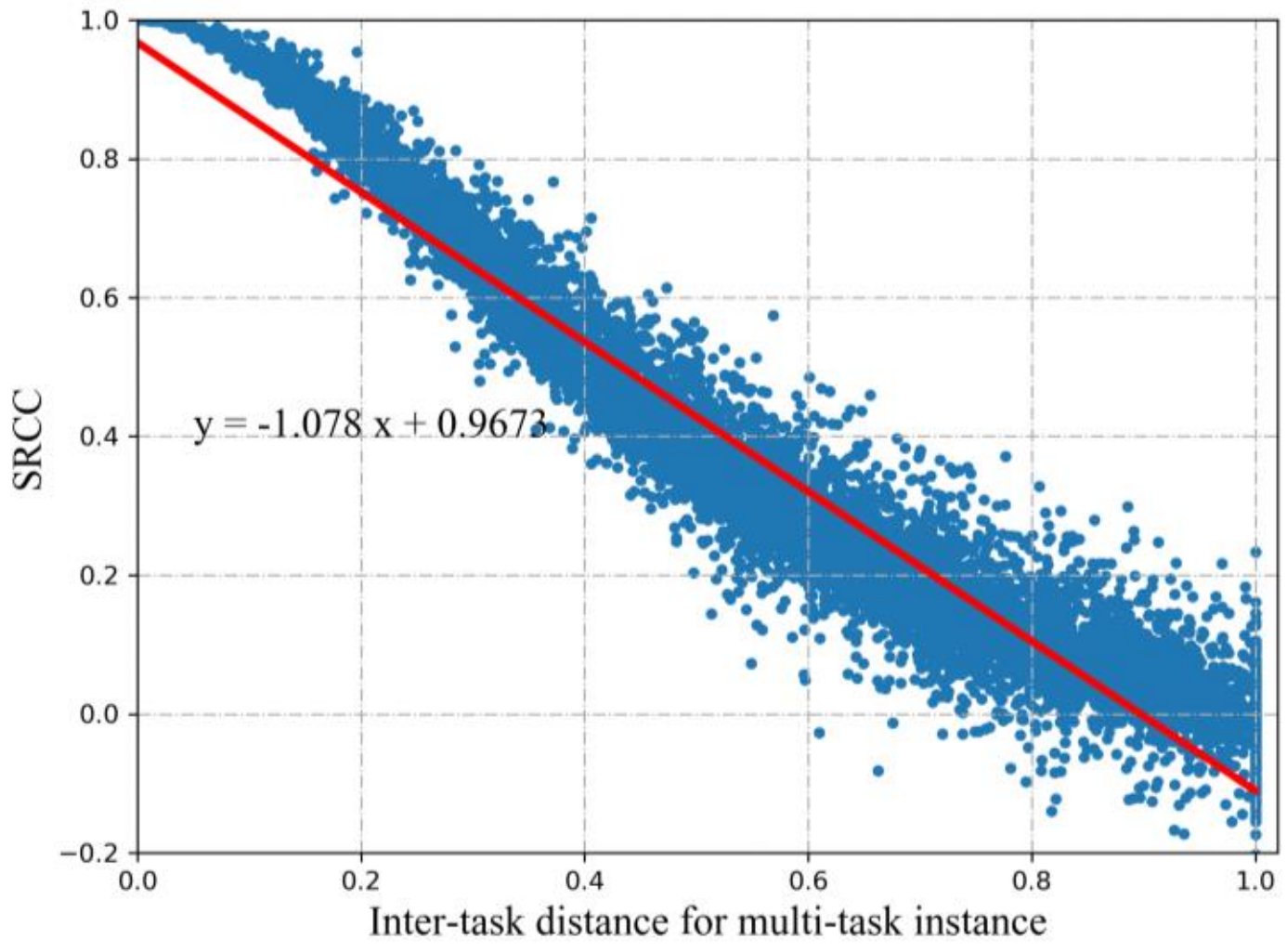

Figure 11: Scatter plot between the inter-task distance measure and Spearman's rank correlation coefficient.

## 9.3 Performance Metric

For comparisons of algorithms, the following performance metrics are considered.

1) Best relative error to $C^*$ (BRE), where, $C^*$ is the optimal or best-so-far value

2) Average relative error to $C^*$ (ARE)

3) Worst relative error to $C^*$ (WRE)

4) Regularized performance score (RPS) (Da et al., 2017)

The $a$ -th algorithm's RPS is

$$RPS_a = \sum_{k=1}^{K} \sum_{l=1}^{L} (I(a,k)_l - \mu_k)/\sigma_k \tag{29}$$

where $I(a,k)_l$ is the best result on the $l$ -th repetition by the $a$ -th algorithm on the $k$ -th problem, $\mu_k$ and $\sigma_k$ are the mean and the standard deviation with respect to the $k$ -th problem over all the repetitions of all algorithms, $L$ is the total repetition of each algorithm on each instance, and $K$ is the number of problems in one instance. The RPS represents the relative results of an algorithm among all algorithms. For our minimization problem, the smaller the RPS, the better the algorithm.

## 9.4 Transferability of High-quality Solutions between Problems

We want to test this hypothesis: The degree to which the high-quality solution to the source problem maintains its optimality in the target problem is highly related to the distance between source and target problems. If it holds when both problems are closed to each other, the high-quality solution of one problem can be directly reused in another problem with a quality guarantee. This could result in dramatic speed-ups in search.

We define the *transferability value* to describe the degree of a performance's non-degradation when the high-quality solution (optimal/best-so-far solution) of Problem 1 (source task) is applied for Problem 2 (target task), and vice versa. To be specific, the transferability value is quantified as the rank of source problem's high quality solution among the target problem's set of solutions.

In the experiment, the source problem's optimal/best-so-far solution, along with the 100,000 randomly generated feasible solutions, are evaluated on target problem through Eqs. (1) – (5). All solutions are sorted in ascending order of the objective values, the makespan. The rank of source problem's optimal/best-so-far solution among the set of solutions for target problem is identified. We compute the transferability value and the inter-task distance for all 12,000 instances. The exact percentage distribution of transferability value is depicted against different groups of distances in Figure 12.

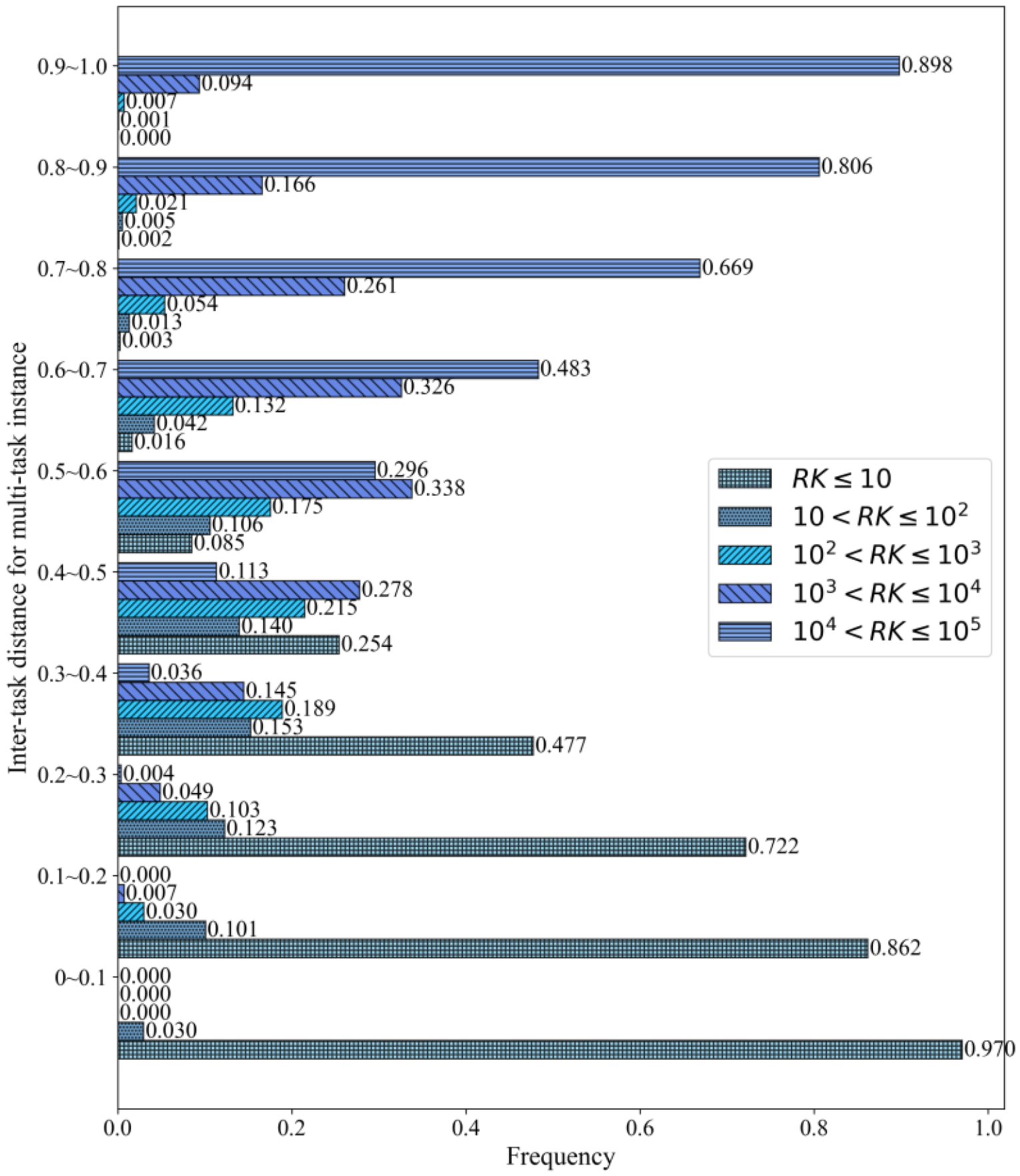

Figure 12: Histogram about the relationship between percentage distribution of transferability value and different groups of distances for 12,000 instances.

As the histogram plot reveals, transferability value has a strong positive dependence on inter-task distance, that is, the closer the distance between the problems, the lower the transferability value (the *RK* in Figure 12). For example, among the instances with the inter-task distance within 0.1 (or 0.2), in nearly 97% (or 86%) instances Problem 1's optimal/best-so-far solution is ranked in the top 10 among the 100,000 solutions for Problem 2. It is an evidence to transfer the source problem's high-quality solutions to the target problem when the inter-task distance is less than 0.2. This would lead to dramatic speed-ups in search.

The performances of optimal/best-so-far solution of Problem 1 are gradually degraded on Problem 2 with the inter-task distance increasing. For example, among the instances with the inter-task distance between 0.9-1.0 (sharing few similarities), there are no instances where Problem 1's optimal/best-so-far solution can be ranked in the top 10 among the 100,000 solutions for Problem 2. Meanwhile, there are nearly 90% instances where Problem 1's optimal/best-so-far solution is ranked in [10,000, 100,000] among the 100,000 solutions for Problem 2. It means the high-quality solution in Problem 1 is severely degenerated in Problem 2 when they are far from each other, and a forcible transfer of complete solution is not beneficial for search.

Furthermore, in investigating the performances of the transfer of high-quality solution (called TranHQ scheme) in depth, we used two high-performance heuristics as comparison algorithms: the NEH (Nawaz, Jr et al. 1983), known as the best constructive heuristic for PFSP, and RandLS, a local search algorithm based on Swap operator with a random initial solution. We selected from Taillard's benchmark 90 instances whose optima are available in OR-library; each is the Problem 1 in a multi-task instance. The Problem 2 is generated by working on the Problem 1 with a replacing probability of 0.05. Each algorithm is independently repeated 50 times on the Problem 2 for each multi-task instance. The makespan and elapsed time by the TranHQ, NEH, and RandLS are exhibited in Table 1. The objective values by the TranHQ and NEH are close and better than those by RandLS. By Wilcoxon signed rank tests, the TranHQ provides an equivalent performance comparing with NEH in statistical sense where p-value for this result is 0.9783. Since, based on elapsed time, TranHQ averages 272 times faster than NEH and 3947 times faster than RandLS, TranHQ is significantly superiority to NEH and RandLS in computation complexity. We can conclude that for very similar problems, high-quality solution transfer would alleviate the time-expensive process for iterative evaluation, and that $10^2 \sim 10^3$ speed up is achieved.

| **Instance Size** | **Inter-task distance** | **Instance Number** | **Objective values (Makespan)** | | | **Elapsed time** | | |
|---|---|---|---|---|---|---|---|---|
| | | | TranHQ | NEH | RandLS | TranHQ | NEH | RandLS |
| 20-5 | 0.10 | 9 | **1240.6** | 1241.6 | 1414.0 | **1.40E-05** | 3.50E-03 | 1.19E-03 |
| 20-10 | 0.13 | 10 | **1561.5** | 1571.8 | 1818.1 | **1.96E-05** | 6.81E-04 | 5.03E-03 |
| 20-20 | 0.12 | 10 | **2295.6** | 2336.8 | 2583.8 | **2.81E-05** | 7.69E-04 | 3.91E-03 |
| 50-5 | 0.16 | 10 | 2821.8 | **2777.4** | 3088.5 | **2.81E-05** | 4.29E-03 | 2.31E-02 |
| 50-10 | 0.14 | 10 | **3106.3** | 3152.5 | 3653.4 | **3.87E-05** | 4.46E-03 | 3.32E-02 |
| 50-20 | 0.14 | 5 | **3842.2** | 3943.2 | 4594.1 | **6.59E-05** | 5.99E-03 | 6.52E-02 |
| 100-5 | 0.16 | 10 | 5370.8 | **5321.1** | 5797.8 | **3.98E-05** | 1.67E-02 | 1.34E-01 |
| 100-10 | 0.14 | 10 | 5798.7 | **5769.0** | 6592.1 | **6.70E-05** | 2.47E-02 | 2.65E-01 |
| 100-20 | 0.13 | 1 | **6504.0** | 6579.0 | 7684.3 | **1.24E-04** | 4.12E-02 | 5.33E-01 |
| 200-10 | 0.15 | 9 | 10954.0 | **10790.8** | 12112.1 | **1.30E-04** | 1.65E-01 | 2.05E+00 |
| 200-20 | 0.14 | 5 | **11764.4** | 11774.6 | 13650.1 | **6.51E-04** | 2.82E-01 | 4.24E+00 |
| 500-20 | 0.14 | 1 | **26717.0** | 26763.0 | 30195.9 | **2.10E-02** | 4.16E+00 | 6.61E+01 |
| **average** | 0.14 | 90 | 4783.9 | **4772.3** | 5401.9 | **3.13E-04** | 8.53E-02 | 1.24E+00 |

Table 1: Comparison between TranHQ, NEH and RandLS in terms of search quality and elapsed time

### 9.5 Comparative Efficiency of the MTCO vs. the Single-task Scatter Search

In comparing the MTCO with the single-task scatter search (STSS), we stripped the multi-task elements (distance measure, problem transformation, and adaptive transfer) from the MTCO to form the STSS. We use the Taillard's instances with 50 jobs and 20 machines (Tai-51, Tai-52, …, Tai-60) since they are particularly difficult to solve (Nowicki and Smutnicki 2006). The replacing probability $p_r$ is ranged from (0, 1] with 10 equal intervals, in effect generating a total 100 instances that also can be used in later experiments. Three parameters were set: size of initial trial solutions $N = 20$, size of reference set $RS = 12$, and maximum number of iteration is $200 \cdot n \cdot (n-1)$. Each instance is independently run 10 times for each algorithm.

As demonstrated in Table 2, MTCO has lower BRE, ARE, WRE and RPS values than STSS. From Figure 13, we see MTCO has a faster convergence speed than STSS, results indicating that the mechanisms of multi-task optimization (distance measure, problem transformation, and adaptive transfer) can indeed provide a bonus to improve performance of the single-task optimization method in both solution quality and searching speed.

| | **Problem 1** | | | **Problem 2** | | | **Instance** | | | |
|---|---|---|---|---|---|---|---|---|---|---|
| | BRE | ARE | WRE | BRE | BRE | ARE | BRE | ARE | WRE | RPS |
| MTCO | 1.63 | 2.17 | 2.64 | 0.45 | 0.99 | 1.47 | 1.04 | 1.58 | 2.06 | -18.47 |
| STSS | 3.46 | 3.94 | 4.30 | 2.33 | 2.76 | 3.12 | 2.89 | 3.35 | 3.71 | 18.47 |

Table 2: Comparison between MTCO and STSS in terms of search quality

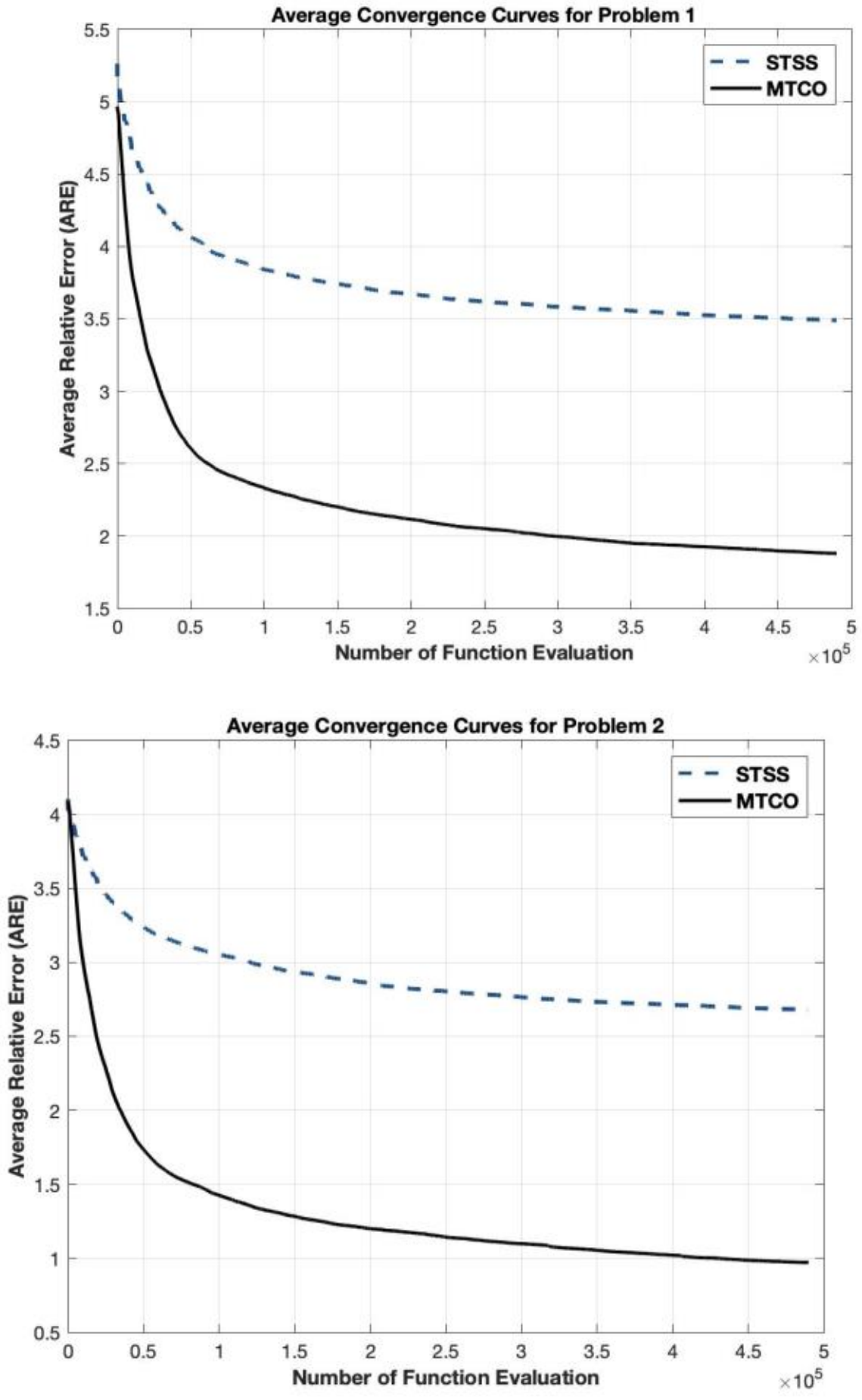

Figure 13: Average convergence curves for Problem 1 and Problem 2.

Furthermore, to quantify the speed-ups of MTCO over STSS in depth, we compared the minimum number of function evaluation for both algorithms to obtain the same high-quality solution. The number of function evaluation by the MTCO and STSS is exhibited in Table 3. Since, based on the number of function evaluation, MTCO averages 3.27 times faster than STSS, MTCO is significantly superiority to STSS in computational cost. For problems with less similarity ( $p_r \geq 0.5$), MTCO averages 3.23 times faster than STSS.

| Replacing probabilities | STSS | MTCO | Speed-up Ratio |
|---|---|---|---|
| 0.1 | 254717.7 | **87421.5** | 2.91 |
| 0.2 | 272428.6 | **73676.1** | 3.70 |
| 0.3 | 267998.5 | **88875.6** | 3.02 |
| 0.4 | 273714.6 | **73811.4** | 3.71 |
| 0.5 | 260167.6 | **82550.3** | 3.15 |
| 0.6 | 265724.1 | **79898.7** | 3.33 |
| 0.7 | 271731.9 | **77644.6** | 3.50 |
| 0.8 | 263378.1 | **91647.7** | 2.87 |
| 0.9 | 268336.2 | **86212.1** | 3.11 |
| 1.0 | 263116.6 | **78544.2** | 3.35 |
| **average** | | | 3.27 |

Table 3: Comparison between MTCO and STSS in terms of number of function evaluation

### 9.6 Effects of the Knowledge Transfer

We investigated the impacts of different types of transferred knowledge on search performance. MTCO-p, MTCO-c and MTCO-e represent the MTCO only transferring one type of knowledge, i.e., partial solutions, complete solution, and solution evolution, respectively. It shows in Table 4 that the MTCO with any type of knowledge transfer performs better than STSS. Among the three types of transfers, MTCO-p is the most effective. By inheriting the invariant knowledge across problems, it reduces the search space by focusing on the non-invariant domain, helpful overall for improved search efficiency.

| | **Problem 1** | | | **Problem 2** | | | **Instance** | | | |
|---|---|---|---|---|---|---|---|---|---|---|
| | BRE | ARE | WRE | BRE | ARE | WRE | BRE | ARE | WRE | RPS |
| MTCO-p | 1.64 | 2.18 | 2.67 | 0.48 | 1.00 | 1.47 | 1.06 | 1.59 | 2.09 | -16.72 |
| MTCO-c | 2.09 | 2.54 | 2.93 | 0.93 | 1.37 | 1.76 | 1.51 | 1.95 | 2.35 | -7.16 |
| MTCO-e | 2.12 | 2.57 | 2.94 | 0.92 | 1.38 | 1.79 | 1.52 | 1.97 | 2.36 | -6.68 |
| STSS | 3.46 | 3.94 | 4.30 | 2.33 | 2.76 | 3.12 | 2.89 | 3.35 | 3.71 | 30.55 |

Table 4: Effects of different types of transferred knowledge on performances

### 9.7 Effects of the Adaptive Transfer Scheme

We investigated the impact of the adaptive transfer scheme guided by the inter-task distance. For comparison, we introduce the MTCO-D0.7 in which all three types of knowledge are transferred in a fixed rate where the inter-task distance $d$ in Eq. (27) is 0.7. This suggests a transfer independent from the actual inter-task distance.

Table 5 shows that on nearly all instances MTCO performs better than MTCO-D0.7. Because MTCO-D0.7 only works well on the instances with their inter-task distances around 0.7, it seems that MTCO-D0.7 cannot promptly adjust the transfer in response to the change of discrepancy between problems. In MTCO, when the actual inter-task similarity is high, the adaptive transfer scheme encourages the transfer to fully use the commonality to foster search. When problems are less related, the adaptive transfer scheme depresses the transfer not to hinder the search. Adjustment of the transfer according to the inter-task distance can improve efficacy.

| **Replacing probabilities** | **Algorithm** | **Problem 1** | | | **Problem 2** | | | **Instance** | | | |
|---|---|---|---|---|---|---|---|---|---|---|---|
| | | BRE | ARE | WRE | BRE | ARE | WRE | BRE | ARE | WRE | RPS |
| 0.1 | MTCO | 1.24 | 1.88 | 2.47 | 0.12 | 0.70 | 1.24 | 0.68 | 1.29 | 1.85 | -7.73 |
| | MTCO-D0.7 | 1.68 | 2.14 | 2.58 | 0.51 | 1.05 | 1.53 | 1.10 | 1.60 | 2.06 | 7.73 |
| 0.2 | MTCO | 1.48 | 2.01 | 2.54 | 0.21 | 0.82 | 1.48 | 0.84 | 1.42 | 2.01 | -4.44 |
| | MTCO-D0.7 | 1.59 | 2.15 | 2.65 | 0.47 | 1.01 | 1.45 | 1.03 | 1.58 | 2.05 | 4.44 |
| 0.3 | MTCO | 1.60 | 2.12 | 2.57 | 0.16 | 0.89 | 1.39 | 0.88 | 1.50 | 1.98 | -4.54 |
| | MTCO-D0.7 | 1.70 | 2.24 | 2.65 | 0.51 | 1.11 | 1.67 | 1.11 | 1.67 | 2.16 | 4.54 |
| 0.4 | MTCO | 1.51 | 2.10 | 2.55 | 0.38 | 1.02 | 1.63 | 0.94 | 1.56 | 2.09 | -5.81 |
| | MTCO-D0.7 | 1.83 | 2.27 | 2.74 | 0.75 | 1.27 | 1.77 | 1.29 | 1.77 | 2.25 | 5.81 |
| 0.5 | MTCO | 1.55 | 2.16 | 2.68 | 0.52 | 0.93 | 1.40 | 1.03 | 1.55 | 2.04 | -2.71 |
| | MTCO-D0.7 | 1.67 | 2.26 | 2.72 | 0.55 | 1.03 | 1.40 | 1.11 | 1.64 | 2.06 | 2.71 |
| 0.6 | MTCO | 1.68 | 2.22 | 2.63 | 0.56 | 1.11 | 1.50 | 1.12 | 1.67 | 2.06 | -2.78 |
| | MTCO-D0.7 | 1.91 | 2.36 | 2.82 | 0.73 | 1.17 | 1.59 | 1.32 | 1.76 | 2.20 | 2.78 |
| 0.7 | MTCO | 1.76 | 2.32 | 2.76 | 0.55 | 1.03 | 1.49 | 1.16 | 1.67 | 2.13 | -1.82 |
| | MTCO-D0.7 | 1.84 | 2.33 | 2.67 | 0.64 | 1.13 | 1.57 | 1.24 | 1.73 | 2.12 | 1.82 |
| 0.8 | MTCO | 1.87 | 2.33 | 2.74 | 0.70 | 1.19 | 1.65 | 1.29 | 1.76 | 2.19 | -2.00 |
| | MTCO-D0.7 | 1.84 | 2.34 | 2.84 | 0.82 | 1.31 | 1.76 | 1.33 | 1.83 | 2.30 | 2.00 |
| 0.9 | MTCO | 1.87 | 2.31 | 2.72 | 0.67 | 1.12 | 1.46 | 1.27 | 1.72 | 2.09 | -1.80 |
| | MTCO-D0.7 | 1.97 | 2.38 | 2.79 | 0.79 | 1.16 | 1.52 | 1.38 | 1.77 | 2.16 | 1.80 |
| 1.0 | MTCO | 1.78 | 2.28 | 2.73 | 0.59 | 1.04 | 1.51 | 1.18 | 1.66 | 2.12 | -1.29 |
| | MTCO-D0.7 | 1.77 | 2.32 | 2.69 | 0.61 | 1.10 | 1.46 | 1.19 | 1.71 | 2.07 | 1.29 |
| **Average** | MTCO | 1.63 | 2.17 | 2.64 | 0.45 | 0.99 | 1.47 | 1.04 | 1.58 | 2.06 | -3.49 |
| | MTCO-D0.7 | 1.78 | 2.28 | 2.72 | 0.64 | 1.13 | 1.57 | 1.21 | 1.71 | 2.14 | 3.49 |

Table 5: Effects of the adaptive transfer scheme.

## 9.8 Effects of Problem Transformation

### 9.8.1 GAIN SIMILARITY OR NOT?

We investigated whether the weakly related problems gain similarity or not after the matching-feature based problem transformation. In Figure 14, the distributions of inter-task distances for instances transformed before and after are depicted against different groups of replacing probabilities $p_r$. It can be revealed from the box plot that when $p_r < 0.5$, there is no significant difference between the original instances' inter-task distances and those of the transformed instances. And when $p_r \geq 0.5$, the inter-task distances of the transformed instances are lower than those of the original instances. As the replacing probability increases (distance between original problems increases), the degree of the change in the distance between the transformed problems also increases. When the distance between the original problems is close to 1 (problems are unrelated), the average distance of the transformed problems is around 0.6.

An interesting observation is that the matching-feature-based problem transformation shows different distance conversion capabilities for problems with different distances. From Figure 14, problems with a distance less than 0.5 cannot be transformed into problems with a smaller distance, and problems with a distance greater than 0.5 cannot be transformed into problems with a distance less than 0.5. It can be concluded that the problem transformation method successfully transforms weakly relevant problems into more relevant ones.

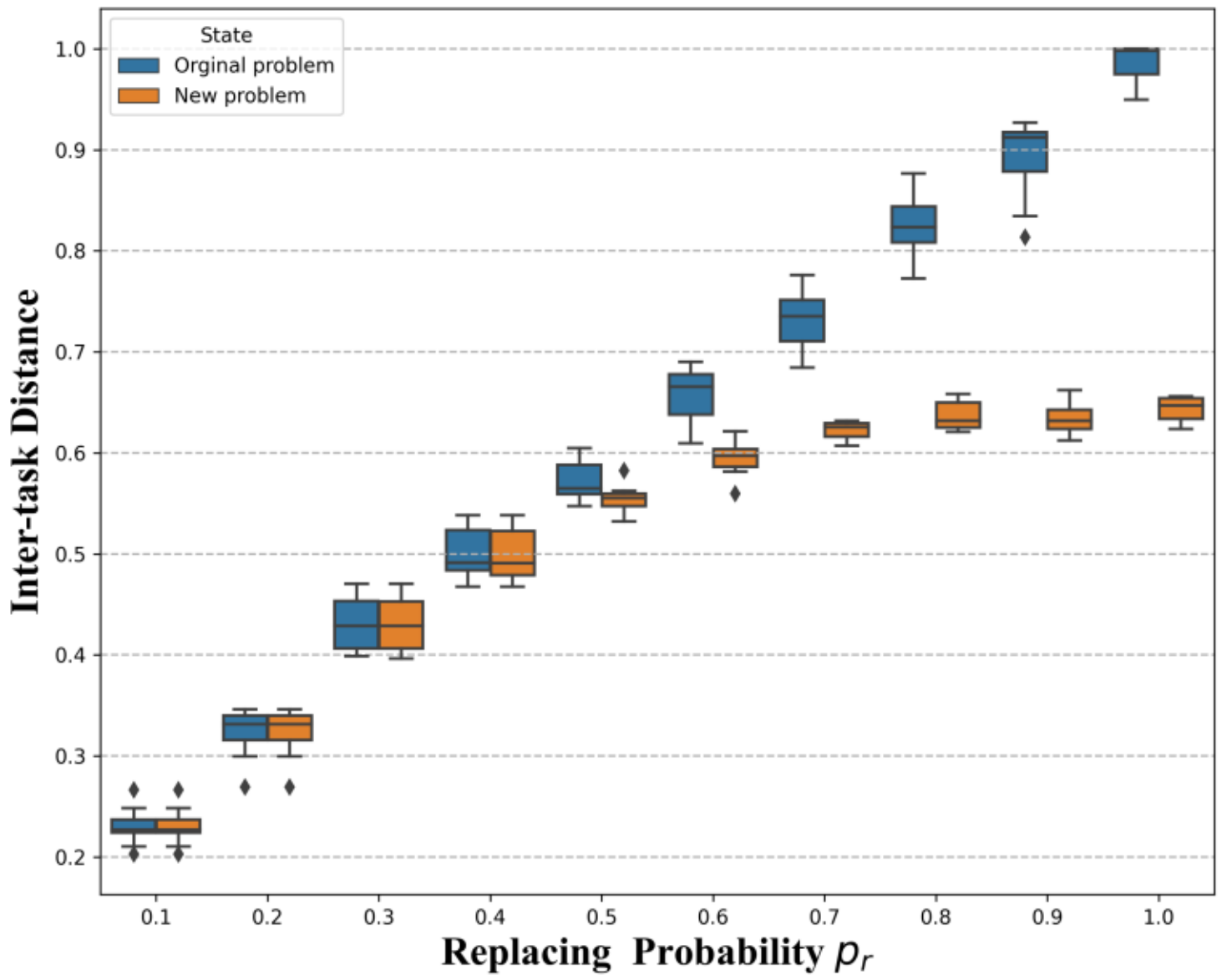

Figure 14: Distributions of inter-task distances for instances transformed before and after are depicted against different groups of replacing probabilities $p_r$.

### 9.8.2 WORKS OR NOT?

We investigated the effect of matching-feature-based problem transformation on search by comparing MTCO with MTCO-noPT (MTCO without problem transformation). Table 6 shows the results when $p_r \leq 0.5$. From the RPS values, it seems that MTCO performs slightly better than MTCO-noPT, but as there is no significant difference between them, the problem transformation cannot significantly improve the search performance when the inter-task distance is less than 0.5. Table 7 provides the results when $p_r > 0.5$. It can be seen that MTCO surpasses MTCO-noPT, suggesting that the problem transformation can significantly improve the search performance by converting the less relevant problems into closer ones. To sum up, mapping irrelevant problems to the latent space where they gain similarity accelerates, through knowledge transfer, the search on the transformed problems.

| | Problem 1 | | | Problem 2 | | | Instance | | | |
|---|---|---|---|---|---|---|---|---|---|---|
| | BRE | ARE | WRE | BRE | ARE | WRE | BRE | ARE | WRE | RPS |
| MTCO | 1.48 | 2.05 | 2.56 | 0.28 | 0.87 | 1.43 | 0.88 | 1.46 | 1.99 | -0.33 |
| MTCO-noPT | 1.53 | 2.06 | 2.56 | 0.35 | 0.90 | 1.40 | 0.94 | 1.48 | 1.98 | 0.33 |

Table 6: Effects of matching-feature-based problem transformation on relevant problems

| | Problem 1 | | | Problem 2 | | | Instance | | | |
|---|---|---|---|---|---|---|---|---|---|---|
| | BRE | ARE | WRE | BRE | ARE | WRE | BRE | ARE | WRE | RPS |
| MTCO | 1.79 | 2.29 | 2.71 | 0.61 | 1.10 | 1.52 | 1.20 | 1.70 | 2.12 | -4.93 |
| MTCO-noPT | 2.00 | 2.45 | 2.84 | 0.81 | 1.25 | 1.63 | 1.40 | 1.85 | 2.23 | 4.93 |

Table 7: Effects of matching-feature-based problem transformation on weakly relevant problems

## 10. Conclusions and Future Work

The quantifiable relationship between scheduling problems is rather unclear, and how to leverage it in combinatorial search remains largely unknown, let alone assessing its effects on search. In tackling such hard questions, this study sought to identify the inter-task distance, and to speed up the search by leveraging this quantifiable commonality. We made several assumptions and see limitations, which should be noted.

*Inter-task distance metric*: It reveals the latent inter-task dependence between flowshop scheduling problems in an analytical and concise manner. Though we validated that it is a faithful measurement for flowshop scheduling problems (traveling salesman problem is a special case), it should be noted that constructing a faithful inter-task distance measure is a challenging topic given the rich variety of representations in combinatorial optimization.

*Multi-modality knowledge*: We mined the common knowledge shared by problems in novel and multifaceted forms, i.e., partial solutions, complete solutions, and solution evolution. These types of common knowledge are problem dependent. It is worthy to further develop more essential mathematical descriptions to reflect the commonality.

*Adaptive transfer scheme*: It determines whether/what to transfer at runtime. We formulated it as an optimal decision problem by observing the inter-task distance and the cumulative votes in favor of (against) transfers. A more comprehensive scheme is required to augment it by the tradeoff between a transfer's benefit and cost, and also by monitoring the search in variations of the landscape when transfer incurs.

*Problem transformation function*: We learned it by a matching-feature-based greedy policy. Though it can succeed in transforming weakly relevant problems into more relevant ones, when it comes to less relevant problems (ones with a distance between [0.7, 1]), it can only shorten their distance to between [0.6, 0.7], not shorter. More advanced matching-feature policy is needed instead of the cheap greedy policy.

*Multi-task combinatorial optimization algorithm*: We formed an effective scatter search based multi-task combinatorial algorithm by incorporating the multi-task techniques, i.e., the explicit inter-task distance metric, the adaptive transfer of multi-modality knowledge, and the matching-feature based problem transformation. It is interesting to validate these multi-task techniques in augmenting the knowledge transfer for various existing algorithms (i.e., case-based reasoning, problem reduction, multi-task Bayesian optimization, evolutionary multitasking, and machine learning).

*Priori knowledge and Boolean learning*: Our multi-task optimization algorithm functions for scheduling problems with different structures and data configurations. However, we disregarded all knowledge once a pair of problems was settled. It is noteworthy that we seek and describe a unifying framework for accumulating the priori knowledge (observations in form of problem structures and performances). To investigate its usage in an accelerated search for new tasks, we feel that Boolean learning with Fourier approximation goes in the direction of an acceptable compromise between exactness and efficiency. Learning with appropriate discrete valued feature maps of limited rank could be a possible research direction.

## Acknowledgements

This research greatly benefited from long discussions with Emeritus Professor Michiel Keyzer (SOW-VU, Faculty of Economics and Business Administration, Vrije Universiteit Amsterdam, The Netherlands) whose many comments were incorporated in the final version. During 2007–2009, co-author B.L. was a junior researcher at SOW-VU under Prof. Keyzer's supervision and gratefully acknowledges his insight. Special thanks go to Professor Ann Marie Ross (University of Chinese Academy of Sciences) for her hard work on polishing the manuscript. B.L. received support from Frontier Science Key Research Program, Chinese Academy of Sciences (QYZDB-SSW-SYS020).

## Appendix A. Nomenclature

| Symbol | Meaning |
|---|---|
| $b$ | real number, $b \in R$ |
| $b^*$ | optimal $b$ |
| $C^*$ | optimal or best-so-far objective value |
| $i$ | job index, $i \in \{1, \cdots, n\}$ |
| $j$ | machine index, $j \in \{1, \cdots, m\}$ |
| $k$ | problem index in one instance |
| $K$ | number of problems in one instance |
| $l$ | repetition index |
| $L$ | total number of repetitions of each algorithm on each instance |
| $m$ | number of machines |
| $n$ | number of jobs |
| $N$ | size of initial trial solutions |
| $p_r$ | replacing probability in generating instance |
| $PT$ | size of the non-invariance set |
| $R$ | real number set |
| $RS$ | size of reference set |
| $t$ | any positive scale value |
| $t^*$ | optimal $t$ |
| $T$ | temperature in annealing process |
| $T_0$ | initial temperature |
| $T_k$ | temperature at generation $k$ |
| $Z$ | sampling size in Spearman's rank correlation coefficient |
| $\lambda$ | cooling coefficient |
| $\mu_k$ | mean with respect to the $k$-th problem over all repetitions of all algorithms |
| $\sigma_k$ | standard deviation with respect to the $k$-th problem over all repetitions of all algorithms |
| $U(1, 99)$ | uniform distribution between 1 and 99 |
| $v$ | voting vector |
| $v(1)$ | cumulative votes in favor of transfers in voting vector $v$ |
| $v(2)$ | cumulative votes against transfers in voting vector $v$ |
| $\pi$ | job permutation, a feasible solution |
| $\pi(i)$ | job in the $i$ -th position of the permutation $\pi$ |
| $d(\pi(i))$ | due date for job $\pi(i)$ |
| $U_i$ | unit penalty for job $\pi(i)$ |
| $C(\pi(i), j)$ | completion time for job $\pi(i)$ on machine $j$ |
| $C(\pi(n), m)$ | maximum completion time for all jobs (i.e., makespan) |
| $\pi_s^{in}, \pi_t^{in}$ | initial solution to the source and target problem |
| $\pi_s^{best}, \pi_t^{best}$ | the best-so-far solution to the source and target problem |
| $\pi_t^{new}$ | the new trial solution to the target problem |
| $X$ | $(n \times n)$ solution matrix |
| $X^T$ | transpose of solution matrix $X$ |
| $X_s^{in}$ | solution matrix for $\pi_s^{in}$ |
| $X_s^{best}$ | solution matrix for $\pi_s^{best}$ |
| $E$ | $(n \times m)$ matrix with all elements equal to one |
| $P$ | $(n \times m)$ problem specification matrix with element $p_{i,j}$ |
| $P_{(i)}$ | the $i$ -th row of $P$ |
| $p_{i,j}$ | element of $P$, processing time of job $i$ on machine $j$ |
| $p_{\pi(i),j}$ | processing time of job $\pi(i)$ on machine $j$ |
| $n_P$ | number of jobs (rows) in the problem specification matrix $P$ |
| $m_P$ | number of machines (columns) in the problem specification matrix $P$ |
| $Q$ | $(n \times m)$ problem specification matrix with element $q_{i,j}$ |
| $Q_{(i)}$ | the $i$ -th row of $Q$ |
| $MF$ | $(n \times n)$ matching-feature matrix with elements of $mf_{i,j}$ |
| $mf_{i,j}$ | Pearson's correlation coefficient between $P_{(i)}$ and $Q_{(j)}$. |
| $S$ | weight matrix with elements $s_{i,j}$ |
| $s_{i,j}$ | Pearson's correlation coefficient |
| $d(f_Q, f_P)$ | distance from $f_Q$ to $f_P$ |
| $\varphi : \pi \to X$ | mapping function from permutation $\pi$ to solution matrix $X$ |
| $\varphi^{-1} : X \to \pi$ | mapping function from solution matrix $X$ to permutation $\pi$ |
| $f$, $f(\pi)$, $f_P$, $f_Q$, $f_P(\pi)$ | objective function |
| $G_p$ | order-isomorphic problems set for function $f_P$ |
| $\Delta E$ | change in the objective function value |
| $H(i)$ | the invariance index for job $i$ |
| $h(i, j)$ | indicative function |
| $O$ | transformation function, $(n \times n)$ permutation matrix |
| $O_P$ | transformation functions for $f_P$ |
| $O^{-1}$ | inverse transformation function |
| $O_\pi$ | solution evolution |
| $I(a,k)_l$ | the best result on the $l$ -th repetition by the $a$ -th algorithm on the $k$ -th problem |

Table A1: Nomenclature

## Appendix B. Proof of Theorem 1

Proof. We first prove $f_{P'}(\pi^0) = t \cdot f_P(\pi^0)$ where $\pi^0 = [1,2,\cdots,n]$ . In other words, $C_{P'}(i, j) = t \cdot C_P(i, j)$ holds for natural numbers, where $C_P(i, j)$ and $C_{P'}(i, j)$ represent the completion time for job $i$ on machine $j$ under problem specification matrix $P$ and $P'$, respectively. Specifically, we show that $C_{P'}(i, j) = t \cdot C_P(i, j)$ holds in the following four successive steps.

(1) Obviously, $C_{P'}(1,1) = p'_{1,1} = t \cdot p_{1,1} = t \cdot C_P(1,1)$ .

(2) Using induction, we prove $C_{P'}(i,1) = t \cdot C_P(i,1)$ holds for all natural numbers as follows.

For the base case, it is required to prove $C_{P'}(2,1) = t \cdot C_P(2,1)$ holds.

According to Eq. (2), we have

$$C_{P'}(2,1) = p'_{2,1} + C_{P'}(1,1) \tag{A.1}$$

then,

$$C_{P'}(2,1) = p'_{2,1} + C_{P'}(1,1) = t \cdot [p_{2,1} + C_P(1,1)] = t \cdot C_P(2,1) \tag{A.2}$$

In the inductive step, assume $C_{P'}(i,1) = t \cdot C_P(i,1)$ holds for natural number $i$ . By the inductive hypothesis, we can get

$$\begin{aligned} C_{P'}(i+1,1) &= p'_{i+1,1} + C_{P'}(i,1) = t \cdot p_{i+1,1} + t \cdot C_P(i,1) \\ &= t \cdot C_P(i+1,1) \end{aligned} \tag{A.3}$$

Thus, once $C_{P'}(i,1) = t \cdot C_P(i,1)$ holds, $C_{P'}(i+1,1) = t \cdot C_P(i+1,1)$ holds for all natural numbers.

(3) It can be obviously proved that $C_{P'}(1, j) = t \cdot C_P(1, j)$ holds for all natural numbers by applying the same procedures in step (2).

(4) We prove $C_{P'}(i, j) = t \cdot C_P(i, j)$ holds for $j = 2, 3, \cdots, m$ and $i = 2, \cdots, n$ as follows.

For the base case, it is required to prove $C_{P'}(l,4-l)=t\cdot C_P(l,4-l)$ is true. From steps (2) and (3), both $C_{P'}(1,3)=t\cdot C_P(1,3)$ and $C_{P'}(3,1)=t\cdot C_P(3,1)$ hold. Then, according to Eq. (4),

$$\begin{aligned} &C_{P'}(2,2)= p'_{2,2}+\max[C_{P'}(1,2),C_{P'}(2,1)]\\ &=t\cdot p_{2,2}+\max[t\cdot C_P(1,2),t\cdot C_P(2,1)]\\ &=t\cdot\{p_{2,2}+\max[C_P(1,2),C_P(2,1)]\}\\ &=t\cdot C_P(2,2)\end{aligned} \tag{A.4}$$

In the inductive step, assume $C_{P'}(l,k-l)=t\cdot C_P(l,k-l)$ holds for some $l$ and $k$ $(2\le l\le k-2)$. We prove $C_{P'}(l,k+1-l)=t\cdot C_P(l,k+1-l)$ holds for $l$ which is a positive integer in the range of $[2,k-1]$.

$$\begin{aligned} &C_{P'}(l,k+1-l)\\ &= p'_{l,k+1-l}+\max[C_{P'}(l-1,k+1-l),C_{P'}(l+1,k-l)]\\ &=t\cdot p_{l,k+1-l}+\max[t\cdot C_P(l-1,k+1-l),t\cdot C_P(l+1,k-l)]\\ &=t\cdot\{p_{l,k+1-l}+\max[C_P(l-1,k+1-l),C_P(l+1,k-l)]\\ &=t\cdot C_P(l,k+1-l)\end{aligned} \tag{A.5}$$

Thus, $C_{P'}(l,k+1-l)=t\cdot C_P(l,k+1-l)$ holds when $C_{P'}(l,k-l)=t\cdot C_P(l,k-l)$ is true. In other words, $C_{P'}(i,j)=t\cdot C_P(i,j)$ is true for $j=2,3,\cdots,m$ and $i=2,\cdots,n$.

Hence, $f_{P'}(\pi^0)=t\cdot f_P(\pi^0)$ holds for $\pi^0=[1,2,\cdots,n]$.

Further, $f_{P'}(\pi)=t\cdot f_P(\pi)$ holds for arbitrary solution $\pi$. We replace the problem matrix $P$ with $X\cdot P$, where $X$ represents the corresponding matrix for solution $\pi$ according to Eq. (6). The proof is identical with the above procedure.

## Appendix C. Proof of Theorem 2

Proof. We prove $C_{P'}(i,j)=C_P(i,j)+(m+n-1)b$ holds where $\pi=[1,2,\cdots,n]$. Specifically, four successive steps are presented to prove Theorem 2.

(1) Obviously, $C_{P'}(1,1)=p'_{1,1}=p_{1,1}+b=C_P(1,1)+b$.

(2) By induction, we prove $C_{P'}(i,1)=C_P(i,1)+i\cdot b$ holds for natural numbers as follows.

For the base case, it is required to prove $C_{P'}(2,1)=C_P(2,1)+2b$ holds.

$$\begin{aligned} C_{P'}(2,1)&= p'_{2,1}+C_{P'}(1,1)=(p_{2,1}+b)+C_P(1,1)+b\\ &=C_P(2,1)+2b\end{aligned} \tag{A.6}$$

In the inductive step, assume $C_{P'}(i,1)=C_P(i,1)+i\cdot b$ holds for $i$. It is needed to prove $C_{P'}(i+1,1)=C_P(i+1,1)+(i+1)\cdot b$ holds. By the inductive hypothesis,

$$\begin{aligned} C_{P'}(i+1,1)&= p'_{i+1,1}+C_{P'}(i,1)=(p_{\mathrm{i}+1,1}+b)+C_P(i,1)+i\cdot b\\ &=C_P(\mathrm{i}+1,1)+b\cdot(i+1)\end{aligned} \tag{A.7}$$

Thus $C_{P'}(i+1,1) = C_P(i+1,1) + (i+1)\cdot b$ is true when $C_{P'}(i,1) = C_P(i,1) + i\cdot b$ holds. That is, $C_{P'}(i,1) = C_P(i,1) + i\cdot b$ holds for natural numbers.

(3) It can be obviously proved that $C_{P'}(1,j) = C_P(1,j) + j\cdot b$ holds for all natural numbers by applying the same procedures in step (2).

(4) We prove $C_{P'}(i,j) = C_P(i,j) + (i+j-1)\cdot b$ holds for $j = 2,3,\cdots,m$ and $i = 1,2,\cdots,n$ as follows.

For the base case, it is required to prove $C_{P'}(l,4-l) = C_P(l,4-l) + 3b$ is true. From steps (2) and (3), both $C_{P'}(1,3) = C_P(1,3) + 3b$ and $C_{P'}(3,1) = t\cdot C_P(3,1) + 3b$ hold. Then, we have,

$$\begin{aligned} C_{P'}(2,2) &= p'_{2,2} + \max[C_{P'}(1,2), C_{P'}(2,1)] \\ &= p_{2,2} + b + \max[C_P(1,2) + 2b, C_P(2,1) + 2b] \\ &= \{p_{2,2} + \max[C_P(1,2), C_P(2,1)]\} + 3b \\ &= C_P(2,2) + 3b \end{aligned} \tag{A.8}$$

In the inductive step, assume $C_{P'}(l,k-l) = C_P(l,k-l) + (k-1)\cdot b$ holds for some $l$ and $k$ $(2 \le l \le k-2)$. We prove $C_{P'}(l,k+1-l) = C_P(l,k+1-l) + k\cdot b$ holds for $l$ which is a positive integer in the range of $[2, k+1-2]$.

$$\begin{aligned} C_{P'}(l,k+1-l) &= p'_{l,k+1-l} + \max[C_{P'}(l-1,k+1-l),\ldots \\ &\quad C_{P'}(l,k-l)] \\ &= p_{l,k+1-l} + b + \max[C_P(l-1,k+1-l)\ldots \\ &\quad + (k-1)\cdot b, C_P(l,k-l) + (k-1)\cdot b] \\ &= \{p_{l,k+1-l} + \max[C_P(l-1,k+1-l),\ldots \\ &\quad C_P(l,k-l)]\} + k\cdot b \\ &= C_P(l,k+1-l) + k\cdot b \end{aligned} \tag{A.9}$$

Thus, $C_{P'}(l,k+1-l) = C_P(l,k+1-l) + k\cdot b$ holds when $C_{P'}(l,k-l) = C_P(l,k-l) + (k-1)\cdot b$ is true. In other words, $C_{P'}(i,j) = C_P(i,j) + (i+j-1)\cdot b$ is true for natural number, suggesting $f_{P'}(\pi) = f_P(\pi) + (m+n-1)\cdot b$ holds.

## Appendix D. Proof of Theorem 3

Proof. Since the transformation matrix $O$ and solution matrix $X$ are arbitrary $n\times n$ permutation matrices, their multiplication $X\cdot O$ is also an $n\times n$ permutation matrix, and it can be mapped to a job permutation. Thus, according to Eq. (7), we can denote $\pi^2 = \varphi^{-1}(X\cdot O)$. The transpose of the permutation matrix $O$ is also a permutation matrix, denoted $O^T$, and we set $\pi^{tran} = \varphi^{-1}(O)$. We denote $\pi^1 = \varphi^{-1}(X)$, and $\pi^0 = [1,2,\cdots,n]$.

Let $P' = O\cdot P$, each element in $P'$ can be calculated as

$$p'_{i,j} = \sum_{g=r} o_{i,r}\, p_{g,j} = p_{\pi^{tran}(i),j} \tag{A.10}$$

According to Eq. (7)

$$\pi^2 = \varphi^{-1}(X \cdot O) = \pi^0 \cdot (X \cdot O)^T = (\pi^0 \cdot O^T) \cdot X^T = \varphi^{-1}(O) \cdot X^T = \pi^{tran} \cdot X^T \tag{A.11}$$

According to Eq. (A.11),

$$p_{\pi^2(i),j} = p_{\pi^{\text{tran}}(\pi^1(i)),j}, \; i = 1,2,\cdots,n; \; j = 1,2,\cdots,m \tag{A.12}$$

According to Eq. (A.10), Eq. (A.12) can be

$$p_{\pi^2(i),j} = p_{\pi^{\text{tran}}(\pi^1(i)),j} = p'_{\pi^1(i),j}, \; i = 1,2,\cdots,n; \; j = 1,2,\cdots,m \tag{A.13}$$

According to Eq. (A.13), it is obvious

$$C_P(\pi^2(1),1) = p_{\pi^2(1),1} = p'_{\pi^1(1),1} = C_{P'}(\pi^1(1),1) \tag{A.14}$$

Using induction, we prove $C_P(\pi^2(i),1) = C_{P'}(\pi^1(i),1)$ holds for all natural numbers as follows.

(1) When $i = 1$, $C_P(\pi^2(i),1) = C_{P'}(\pi^1(i),1)$ is true according to Eq. (A.14).

(2) When $i > 1$, assume $C_P(\pi^2(i),1) = C_{P'}(\pi^1(i),1)$ holds for natural number $i$. By the inductive hypothesis, we can get

$$C_P(\pi^2(i+1),1) = p_{\pi^2(i+1),1} + C_P(\pi^2(i),1) = p'_{\pi^1(i+1),1} + C_{P'}(\pi^1(i),1) = C_{P'}(\pi^1(i+1),1) \tag{A.15}$$

Thus, once $C_P(\pi^2(i),1) = C_{P'}(\pi^1(i),1)$ holds, $C_P(\pi^2(i+1),1) = C_{P'}(\pi^1(i+1),1)$ holds for all natural numbers.

(3) It can be obviously proved that

$$C_P(\pi^2(1),j) = C_{P'}(\pi^1(1),j), \quad j = 2,3,\cdots,m \tag{A.16}$$

holds for all natural numbers by applying the same procedures in step (2).

According to Eq. (A.15) and Eq. (A.16), by induction it can be proved

$$C_P(\pi^2(i),j) = C_{P'}(\pi^1(i),j) \quad i = 2,3,\cdots,n; \; j = 2,3,\cdots,m \tag{A.17}$$

When $i = n$ and $j = m$, Eq. (A.17) is

$$f_P(\pi^2) = C_P(\pi^2(n),m) = C_{P'}(\pi^1(n),m) = f_{P'}(\pi^1) \tag{A.18}$$

Thus, when $P' = O \cdot P$, $f_P(\pi^2) = f_{P'}(\pi^1)$ holds, and consequently $f_{P'}(\varphi^{-1}(X)) = f_P(\varphi^{-1}(X \cdot O))$ holds for any solution matrix.

## Appendix E. Proof of Theorem 4

Proof. Since the transformation function $O$ is an arbitrary $n \times n$ permutation matrix, its inverse matrix $O^{-1}$ equals to its transposed matrix $O^T$, and remains a permutation matrix. Here, to get Theorem 4, we just need to replace the permutation matrix $O$ with $O^{-1}$, and interchange $P$ and $P'$ in theorem

3. Like the proof in Theorem 3, Theorem 4 can be proved in the same way. However, we don't repeat the proof process here.

## Appendix F. Normalization of the Inter-task Distance

The $Q^*$ in Eq. (19) and $P^*$ in Eq. (20) could be represented as two points in a space with dimension of $n \cdot m$ (see Figure A1). According to Eq. (18), in the aspect of geometric meaning, the inter-task distance could be interpreted as the minimum distance from point $Q^*$ to the ray $tP^* (t \geq 0)$. $\theta$ is the angle between the ray $tP^* (t \geq 0)$ and ray $t \cdot Q^* (t \geq 0)$.

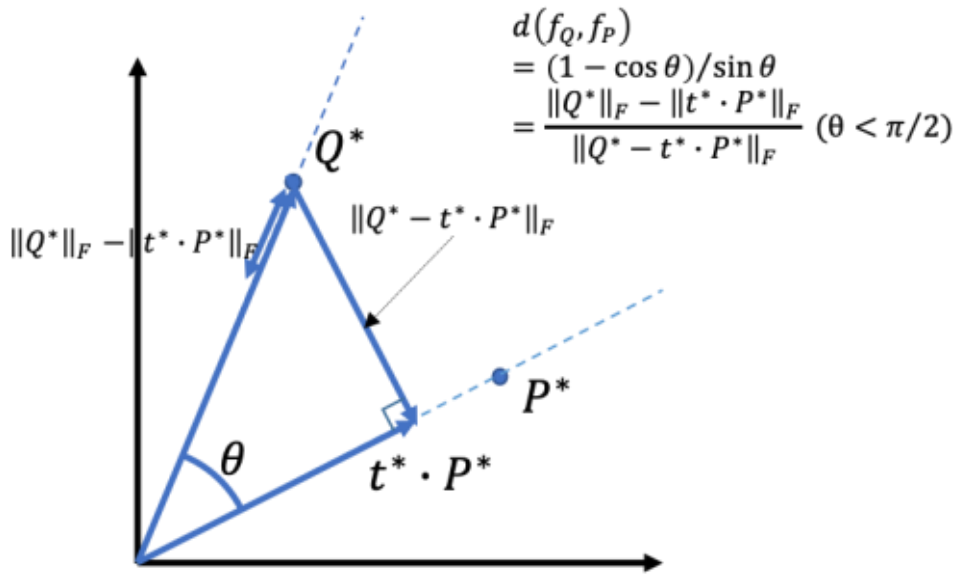

Figure A1: Geometric interpretation: the inter-task distance when angle is less than $\pi/2$.

When $\theta < \pi/2$ (see Figure A1), the normalized distance Eq. (22) can be represented as:

$$d(f_Q, f_P) = (1 - \cos\theta)/\sin\theta \tag{A.19}$$

where,

$$\sin\theta = \left\|Q^* - t^* \cdot P^*\right\|_F \Big/ \left\|Q^*\right\|_F \tag{A.20}$$

$$\cos\theta = t^* \cdot \left\|P^*\right\|_F \Big/ \left\|Q^*\right\|_F \tag{A.21}$$

Obviously, in the interval of $[0, \pi/2]$, the function $(1 - \cos\theta)/\sin\theta$ is increasing with $\theta$ and the inter-task distance Eq. (22) increases from 0 to 1. Distance value of 1 suggests completely unrelated PFSPs (i.e., PFSPs shares no similarities), while value of 0 suggests identical PFSPs, i.e., there exists order-isomorphic relationship between PFSPs according to Definition 1.

When $\theta \geq \pi/2$ (see Figure A2), the minimum point in the ray $tP^* (t \geq 0)$ to the point $Q^*$ is the origin, i.e., $t^* = 0$. In this regard, the inter-task distance $\left\|Q^* - t^* \cdot P^*\right\|_F$ is equal to $\left\|Q^*\right\|_F$, and $d(f_Q, f_P) = 1$ according to Eq. (A.19).

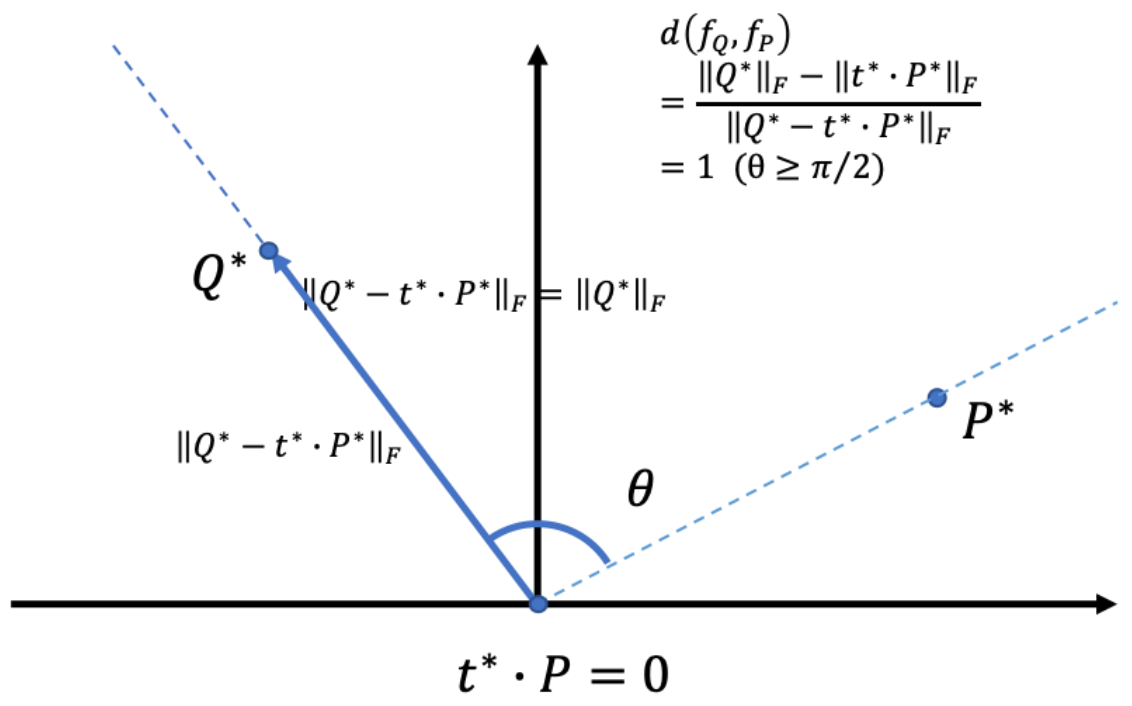

Figure A2: Geometric interpretation: the inter-task distance when angle is greater than or equal to $\pi/2$.

## Appendix G. Computational Complexity of Inter-task Distance Metric

The computational complexity is estimated by the bounds on the number of floating-point operations. Table A2 exhibits the number of addition and multiplication operations for each component in calculating the inter-task distance. There are $(10 \cdot n \cdot m - 3)$ additions and $(6 \cdot n \cdot m + 10)$ multiplications. The calculations of $\|Q^*\|_F$, $\|P^*\|_F$ and $\|Q^* - t^* \cdot P^*\|_F$ involve three root square operations. The inter-task distance can be computed in $O(m \cdot n)$.

| Components | Number of Addition Operations | Number of Multiplication Operations |
|---|---|---|
| $\sum_{i=1}^{n}\sum_{j=1}^{n} p_{i,j}$ | $nm-1$ | 0 |
| $\sum_{i=1}^{n}\sum_{j=1}^{n} q_{i,j}$ | $nm-1$ | 0 |
| $\sum_{i=1}^{n}\sum_{j=1}^{n} q_{i,j} p_{i,j}$ | $nm-1$ | $nm$ |
| $\sum_{i=1}^{n}\sum_{j=1}^{n} (p_{i,j})^2$ | $nm-1$ | $nm$ |
| $n \cdot m$ | 0 | 1 |
| $P^* = P - \frac{1}{n \cdot m} \cdot \sum_{j=1}^{n}\sum_{i=1}^{m} p_{i,j} \cdot E$ | $nm$ | 1 |
| $Q^* = Q - \frac{1}{n \cdot m} \cdot \sum_{j=1}^{n}\sum_{i=1}^{m} q_{i,j} \cdot E$ | $nm$ | 1 |
| $n \cdot m \cdot \sum_{j=1}^{n}\sum_{i=1}^{m} q_{i,j} \cdot p_{i,j} - (\sum_{j=1}^{n}\sum_{i=1}^{m} q_{i,j}) \cdot (\sum_{j=1}^{n}\sum_{i=1}^{m} p_{i,j})$ | 1 | 2 |
| $t^*$ | 2 | 3 |
| $\left\|Q^* - t^* \cdot P^*\right\|_F$ | $2nm-1$ | $2nm$ |
| $\left\|Q^*\right\|_F$ | $nm-1$ | $nm$ |
| $\left\|P^*\right\|_F$ | $nm-1$ | $nm$ |
| $d(f_Q, f_P) = \frac{\left\|Q^*\right\|_F - t^* \cdot \left\|P^*\right\|_F}{\left\|Q^* - t^* \cdot P^*\right\|_F}$ | 1 | 2 |
| **Sum** | $10nm-3$ | $6nm+10$ |

Table A2: Computational Complexity of Inter-task Distance

## Appendix H. An Illustrative Example for Computing the Invariance Index

In this illustrative example, $\pi^p = [3,4,2,5,6,1]$ and $\pi^q = [1,4,2,6,5,3]$. We start from job 1. It is observed that job 1 succeeds job 2 in $\pi^p$, while job 1 precedes job 2 in $\pi^q$. Thus, $h(1,2) = 0$. Following the same way, the invariance indexes for all jobs are obtained in Table A3.

| **Job** | 1 | 2 | 3 | 4 | 5 | 6 |
|---|---|---|---|---|---|---|
| 1 | - | 0 | 0 | 0 | 0 | 0 |
| 2 | 0 | - | 0 | 1 | 1 | 1 |
| 3 | 0 | 0 | - | 0 | 0 | 0 |
| 4 | 0 | 1 | 0 | - | 1 | 1 |
| 5 | 0 | 1 | 0 | 1 | - | 0 |
| 6 | 0 | 1 | 0 | 1 | 0 | - |
| **H** | 0 | 0.6 | 0 | 0.6 | 0.4 | 0.4 |

Table A3: An Illustrative Example for Computing the Invariance Index

## Appendix I. Precedence-based Distance Metric to Measure the Distance between Optima

The distance between optima for two problems is calculated by the well-known precedence-based distance metric (Reeves & Yamada, 1998). It counts the minimum number of interchange moves, necessary for transforming $\pi_1$ to $\pi_2$ as follows:

$$count(\pi_1, \pi_2) = |\{i, j\} : [\pi_1^{-1}(i) - \pi_1^{-1}(j)] \cdot [\pi_2^{-1}(i) - \pi_2^{-1}(j)] < 0, \ 1 \le i < j \le n| \tag{A.22}$$

where, $count(\pi_1, \pi_2)$ counts the number of job pairs satisfying the condition in Eq. (A.22), $\pi_1^{-1}$ and $\pi_2^{-1}$ are the jobs' precedence ranks in the job sequences $\pi_1$ and $\pi_2$, respectively, and the element in $\pi^{-1}(i)$ denotes the position of job $i$ in solution $\pi$.

Further, $count(\pi_1, \pi_2)$ is divided by $n(n-1)/2$ to obtain a normalized precedence-based distance. $n(n-1)/2$ is the maximum number of interchange for transform between two solutions,

$$D(\pi_1, \pi_2) = 2 \cdot count(\pi_1, \pi_2)/(n(n-1)) \tag{A.23}$$

## Appendix J. Computational Results of Inter-task Distance for All Instances

By pairing instances with identical dimensions from Taillard's benchmark, multiple instances are generated and their inter-task distances are measured. The results for all instances with $[20, 50, 100, 200, 500]$ jobs and $[5, 10, 20]$ machines are shown in Tables A.4-A.14. Distance value one suggests unrelated PFSPs. Almost all instances from the single-task Taillard's benchmark are unrelated or faintly related in terms of the inter-task distance metric.

| | Tai1 | Tai2 | Tai3 | Tai4 | Tai5 | Tai6 | Tai7 | Tai8 | Tai9 | Tai10 |
|---|---|---|---|---|---|---|---|---|---|---|
| Tai1 | 0.00 | 1.00 | 1.00 | 0.91 | 1.00 | 0.94 | 0.93 | 1.00 | 1.00 | 0.90 |
| Tai2 | 1.00 | 0.00 | 1.00 | 0.86 | 1.00 | 0.99 | 0.89 | 1.00 | 1.00 | 1.00 |
| Tai3 | 1.00 | 1.00 | 0.00 | 1.00 | 0.99 | 0.91 | 0.98 | 1.00 | 0.96 | 1.00 |
| Tai4 | 0.91 | 0.86 | 1.00 | 0.00 | 1.00 | 0.86 | 0.85 | 1.00 | 1.00 | 1.00 |
| Tai5 | 1.00 | 1.00 | 0.99 | 1.00 | 0.00 | 0.99 | 1.00 | 0.86 | 0.91 | 1.00 |
| Tai6 | 0.94 | 0.99 | 0.91 | 0.86 | 0.99 | 0.00 | 0.83 | 1.00 | 0.96 | 1.00 |
| Tai7 | 0.93 | 0.89 | 0.98 | 0.85 | 1.00 | 0.83 | 0.00 | 1.00 | 1.00 | 0.94 |
| Tai8 | 1.00 | 1.00 | 1.00 | 1.00 | 0.86 | 1.00 | 1.00 | 0.00 | 0.99 | 1.00 |
| Tai9 | 1.00 | 1.00 | 0.96 | 1.00 | 0.91 | 0.96 | 1.00 | 0.99 | 0.00 | 0.95 |
| Tai10 | 0.90 | 1.00 | 1.00 | 1.00 | 1.00 | 1.00 | 0.94 | 1.00 | 0.95 | 0.00 |

Table A.4 The Inter-task Distance between Taillard's Instances with 20 Jobs and 5 Machines.

| | Tai11 | Tai12 | Tai13 | Tai14 | Tai15 | Tai16 | Tai17 | Tai18 | Tai19 | Tai20 |
|---|---|---|---|---|---|---|---|---|---|---|
| Tai11 | 0.00 | 1.00 | 0.98 | 1.00 | 0.96 | 1.00 | 0.98 | 0.99 | 0.93 | 0.89 |
| Tai12 | 1.00 | 0.00 | 1.00 | 0.93 | 1.00 | 1.00 | 0.97 | 1.00 | 0.91 | 1.00 |
| Tai13 | 0.98 | 1.00 | 0.00 | 0.99 | 1.00 | 1.00 | 0.85 | 0.90 | 1.00 | 1.00 |
| Tai14 | 1.00 | 0.93 | 0.99 | 0.00 | 1.00 | 0.93 | 0.95 | 1.00 | 0.99 | 1.00 |
| Tai15 | 0.96 | 1.00 | 1.00 | 1.00 | 0.00 | 1.00 | 1.00 | 1.00 | 0.89 | 1.00 |
| Tai16 | 1.00 | 1.00 | 1.00 | 0.93 | 1.00 | 0.00 | 1.00 | 0.91 | 1.00 | 1.00 |
| Tai17 | 0.98 | 0.97 | 0.85 | 0.95 | 1.00 | 1.00 | 0.00 | 0.88 | 0.92 | 0.99 |
| Tai18 | 0.99 | 1.00 | 0.90 | 1.00 | 1.00 | 0.91 | 0.88 | 0.00 | 1.00 | 1.00 |
| Tai19 | 0.93 | 0.91 | 1.00 | 0.99 | 0.89 | 1.00 | 0.92 | 1.00 | 0.00 | 0.86 |
| Tai20 | 0.89 | 1.00 | 1.00 | 1.00 | 1.00 | 1.00 | 0.99 | 1.00 | 0.86 | 0.00 |

Table A.5 The Inter-task Distance between Taillard’s Instances with 20 Jobs and 10 Machines.

| | Tai21 | Tai22 | Tai23 | Tai24 | Tai25 | Tai26 | Tai27 | Tai28 | Tai29 | Tai30 |
|---|---|---|---|---|---|---|---|---|---|---|
| Tai21 | 0.00 | 0.91 | 1.00 | 0.97 | 0.93 | 0.95 | 0.91 | 1.00 | 1.00 | 1.00 |
| Tai22 | 0.91 | 0.00 | 1.00 | 1.00 | 1.00 | 1.00 | 1.00 | 1.00 | 0.95 | 0.97 |
| Tai23 | 1.00 | 1.00 | 0.00 | 1.00 | 1.00 | 1.00 | 1.00 | 1.00 | 1.00 | 1.00 |
| Tai24 | 0.97 | 1.00 | 1.00 | 0.00 | 1.00 | 0.97 | 0.97 | 1.00 | 1.00 | 1.00 |
| Tai25 | 0.93 | 1.00 | 1.00 | 1.00 | 0.00 | 0.98 | 0.91 | 1.00 | 1.00 | 1.00 |
| Tai26 | 0.95 | 1.00 | 1.00 | 0.97 | 0.98 | 0.00 | 0.91 | 0.95 | 0.99 | 0.96 |
| Tai27 | 0.91 | 1.00 | 1.00 | 0.97 | 0.91 | 0.91 | 0.00 | 0.98 | 0.98 | 0.95 |
| Tai28 | 1.00 | 1.00 | 1.00 | 1.00 | 1.00 | 0.95 | 0.98 | 0.00 | 0.95 | 1.00 |
| Tai29 | 1.00 | 0.95 | 1.00 | 1.00 | 1.00 | 0.99 | 0.98 | 0.95 | 0.00 | 1.00 |
| Tai30 | 1.00 | 0.97 | 1.00 | 1.00 | 1.00 | 0.96 | 0.95 | 1.00 | 1.00 | 0.00 |

Table A.6 The Inter-task Distance between Taillard’s Instances with 20 Jobs and 20 Machines.

| | Tai31 | Tai32 | Tai33 | Tai34 | Tai35 | Tai36 | Tai37 | Tai38 | Tai39 | Tai40 |
|---|---|---|---|---|---|---|---|---|---|---|
| Tai31 | 0.00 | 1.00 | 1.00 | 1.00 | 1.00 | 1.00 | 1.00 | 0.99 | 1.00 | 1.00 |
| Tai32 | 1.00 | 0.00 | 0.97 | 0.98 | 1.00 | 1.00 | 0.87 | 0.89 | 0.99 | 1.00 |
| Tai33 | 1.00 | 0.97 | 0.00 | 0.95 | 1.00 | 0.98 | 0.99 | 0.98 | 0.89 | 1.00 |
| Tai34 | 1.00 | 0.98 | 0.95 | 0.00 | 0.93 | 0.94 | 1.00 | 0.98 | 1.00 | 0.94 |
| Tai35 | 1.00 | 1.00 | 1.00 | 0.93 | 0.00 | 1.00 | 0.99 | 0.95 | 1.00 | 0.98 |
| Tai36 | 1.00 | 1.00 | 0.98 | 0.94 | 1.00 | 0.00 | 0.95 | 1.00 | 1.00 | 0.98 |
| Tai37 | 1.00 | 0.87 | 0.99 | 1.00 | 0.99 | 0.95 | 0.00 | 0.92 | 1.00 | 1.00 |
| Tai38 | 0.99 | 0.89 | 0.98 | 0.98 | 0.95 | 1.00 | 0.92 | 0.00 | 1.00 | 1.00 |
| Tai39 | 1.00 | 0.99 | 0.89 | 1.00 | 1.00 | 1.00 | 1.00 | 1.00 | 0.00 | 1.00 |
| Tai40 | 1.00 | 1.00 | 1.00 | 0.94 | 0.98 | 0.98 | 1.00 | 1.00 | 1.00 | 0.00 |

Table A.7 The Inter-task Distance between Taillard’s Instances with 50 Jobs and 5 Machines.

| | Tai41 | Tai42 | Tai43 | Tai44 | Tai45 | Tai46 | Tai47 | Tai48 | Tai49 | Tai50 |
|---|---|---|---|---|---|---|---|---|---|---|
| Tai41 | 0.00 | 0.97 | 0.99 | 0.99 | 0.91 | 1.00 | 0.93 | 1.00 | 1.00 | 0.98 |
| Tai42 | 0.97 | 0.00 | 1.00 | 0.95 | 0.99 | 0.96 | 1.00 | 1.00 | 0.98 | 1.00 |
| Tai43 | 0.99 | 1.00 | 0.00 | 1.00 | 1.00 | 0.97 | 1.00 | 0.99 | 1.00 | 1.00 |
| Tai44 | 0.99 | 0.95 | 1.00 | 0.00 | 0.97 | 0.94 | 1.00 | 1.00 | 0.96 | 0.96 |
| Tai45 | 0.91 | 0.99 | 1.00 | 0.97 | 0.00 | 1.00 | 0.98 | 1.00 | 0.97 | 1.00 |
| Tai46 | 1.00 | 0.96 | 0.97 | 0.94 | 1.00 | 0.00 | 1.00 | 1.00 | 0.97 | 0.99 |
| Tai47 | 0.93 | 1.00 | 1.00 | 1.00 | 0.98 | 1.00 | 0.00 | 0.94 | 1.00 | 1.00 |
| Tai48 | 1.00 | 1.00 | 0.99 | 1.00 | 1.00 | 1.00 | 0.94 | 0.00 | 1.00 | 1.00 |
| Tai49 | 1.00 | 0.98 | 1.00 | 0.96 | 0.97 | 0.97 | 1.00 | 1.00 | 0.00 | 1.00 |
| Tai50 | 0.98 | 1.00 | 1.00 | 0.96 | 1.00 | 0.99 | 1.00 | 1.00 | 1.00 | 0.00 |

Table A.8 The Inter-task Distance between Taillard's Instances with 50 Jobs and 10 Machines.

| | Tai61 | Tai62 | Tai63 | Tai64 | Tai65 | Tai66 | Tai67 | Tai68 | Tai69 | Tai70 |
|---|---|---|---|---|---|---|---|---|---|---|
| Tai61 | 0.00 | 0.98 | 1.00 | 0.92 | 1.00 | 0.98 | 0.98 | 1.00 | 0.99 | 0.94 |
| Tai62 | 0.98 | 0.00 | 0.93 | 0.98 | 1.00 | 1.00 | 0.93 | 1.00 | 1.00 | 1.00 |
| Tai63 | 1.00 | 0.93 | 0.00 | 0.99 | 1.00 | 1.00 | 0.97 | 0.97 | 1.00 | 1.00 |
| Tai64 | 0.92 | 0.98 | 0.99 | 0.00 | 0.98 | 1.00 | 1.00 | 1.00 | 1.00 | 0.95 |
| Tai65 | 1.00 | 1.00 | 1.00 | 0.98 | 0.00 | 0.95 | 1.00 | 1.00 | 1.00 | 0.98 |
| Tai66 | 0.98 | 1.00 | 1.00 | 1.00 | 0.95 | 0.00 | 1.00 | 0.94 | 0.96 | 0.98 |
| Tai67 | 0.98 | 0.93 | 0.97 | 1.00 | 1.00 | 1.00 | 0.00 | 1.00 | 1.00 | 0.97 |
| Tai68 | 1.00 | 1.00 | 0.97 | 1.00 | 1.00 | 0.94 | 1.00 | 0.00 | 0.92 | 1.00 |
| Tai69 | 0.99 | 1.00 | 1.00 | 1.00 | 1.00 | 0.96 | 1.00 | 0.92 | 0.00 | 1.00 |
| Tai70 | 0.94 | 1.00 | 1.00 | 0.95 | 0.98 | 0.98 | 0.97 | 1.00 | 1.00 | 0.00 |

Table A.9 The Inter-task Distance between Taillard's Instances with 100 Jobs and 5 Machines.

| | Tai71 | Tai72 | Tai73 | Tai74 | Tai75 | Tai76 | Tai77 | Tai78 | Tai79 | Tai80 |
|---|---|---|---|---|---|---|---|---|---|---|
| Tai71 | 0.00 | 1.00 | 0.96 | 1.00 | 1.00 | 0.98 | 0.99 | 1.00 | 1.00 | 1.00 |
| Tai72 | 1.00 | 0.00 | 1.00 | 1.00 | 1.00 | 0.97 | 1.00 | 1.00 | 1.00 | 0.97 |
| Tai73 | 0.96 | 1.00 | 0.00 | 0.98 | 1.00 | 1.00 | 0.97 | 1.00 | 0.97 | 1.00 |
| Tai74 | 1.00 | 1.00 | 0.98 | 0.00 | 0.94 | 1.00 | 0.96 | 1.00 | 1.00 | 0.98 |
| Tai75 | 1.00 | 1.00 | 1.00 | 0.94 | 0.00 | 1.00 | 1.00 | 1.00 | 1.00 | 0.96 |
| Tai76 | 0.98 | 0.97 | 1.00 | 1.00 | 1.00 | 0.00 | 0.97 | 0.97 | 1.00 | 1.00 |
| Tai77 | 0.99 | 1.00 | 0.97 | 0.96 | 1.00 | 0.97 | 0.00 | 0.98 | 1.00 | 0.96 |
| Tai78 | 1.00 | 1.00 | 1.00 | 1.00 | 1.00 | 0.97 | 0.98 | 0.00 | 1.00 | 1.00 |
| Tai79 | 1.00 | 1.00 | 0.97 | 1.00 | 1.00 | 1.00 | 1.00 | 1.00 | 0.00 | 1.00 |
| Tai80 | 1.00 | 0.97 | 1.00 | 0.98 | 0.96 | 1.00 | 0.96 | 1.00 | 1.00 | 0.00 |

Table A.10 The Inter-task Distance between Taillard's Instances with 100 Jobs and 10 Machines.

| | Tai81 | Tai82 | Tai83 | Tai84 | Tai85 | Tai86 | Tai87 | Tai88 | Tai89 | Tai90 |
|---|---|---|---|---|---|---|---|---|---|---|
| Tai81 | 0.00 | 1.00 | 0.99 | 1.00 | 1.00 | 0.99 | 1.00 | 0.99 | 1.00 | 1.00 |
| Tai82 | 1.00 | 0.00 | 1.00 | 1.00 | 0.99 | 0.97 | 1.00 | 0.97 | 1.00 | 1.00 |
| Tai83 | 0.99 | 1.00 | 0.00 | 1.00 | 1.00 | 1.00 | 0.97 | 1.00 | 0.96 | 0.98 |
| Tai84 | 1.00 | 1.00 | 1.00 | 0.00 | 0.97 | 1.00 | 1.00 | 1.00 | 1.00 | 1.00 |
| Tai85 | 1.00 | 0.99 | 1.00 | 0.97 | 0.00 | 0.99 | 1.00 | 0.99 | 1.00 | 1.00 |
| Tai86 | 0.99 | 0.97 | 1.00 | 1.00 | 0.99 | 0.00 | 1.00 | 0.98 | 1.00 | 0.98 |
| Tai87 | 1.00 | 1.00 | 0.97 | 1.00 | 1.00 | 1.00 | 0.00 | 0.98 | 0.99 | 1.00 |
| Tai88 | 0.99 | 0.97 | 1.00 | 1.00 | 0.99 | 0.98 | 0.98 | 0.00 | 1.00 | 1.00 |
| Tai89 | 1.00 | 1.00 | 0.96 | 1.00 | 1.00 | 1.00 | 0.99 | 1.00 | 0.00 | 1.00 |
| Tai90 | 1.00 | 1.00 | 0.98 | 1.00 | 1.00 | 0.98 | 1.00 | 1.00 | 1.00 | 0.00 |

Table A.11 The Inter-task Distance between Taillard's Instances with 100 Jobs and 20 Machines.

| | Tai91 | Tai92 | Tai93 | Tai94 | Tai95 | Tai96 | Tai97 | Tai98 | Tai99 | Tai100 |
|---|---|---|---|---|---|---|---|---|---|---|
| Tai91 | 0.00 | 1.00 | 0.98 | 1.00 | 0.99 | 0.98 | 1.00 | 1.00 | 0.98 | 1.00 |
| Tai92 | 1.00 | 0.00 | 1.00 | 0.98 | 1.00 | 1.00 | 0.98 | 1.00 | 0.98 | 0.97 |
| Tai93 | 0.98 | 1.00 | 0.00 | 1.00 | 1.00 | 0.97 | 0.98 | 1.00 | 1.00 | 0.99 |
| Tai94 | 1.00 | 0.98 | 1.00 | 0.00 | 0.97 | 0.99 | 1.00 | 0.99 | 1.00 | 1.00 |
| Tai95 | 0.99 | 1.00 | 1.00 | 0.97 | 0.00 | 0.99 | 1.00 | 0.98 | 0.98 | 0.98 |
| Tai96 | 0.98 | 1.00 | 0.97 | 0.99 | 0.99 | 0.00 | 1.00 | 1.00 | 1.00 | 1.00 |
| Tai97 | 1.00 | 0.98 | 0.98 | 1.00 | 1.00 | 1.00 | 0.00 | 1.00 | 1.00 | 1.00 |
| Tai98 | 1.00 | 1.00 | 1.00 | 0.99 | 0.98 | 1.00 | 1.00 | 0.00 | 0.99 | 1.00 |
| Tai99 | 0.98 | 0.98 | 1.00 | 1.00 | 0.98 | 1.00 | 1.00 | 0.99 | 0.00 | 0.98 |
| Tai100 | 1.00 | 0.97 | 0.99 | 1.00 | 0.98 | 1.00 | 1.00 | 1.00 | 0.98 | 0.00 |

Table A.12 The Inter-task Distance between Taillard's Instances with 200 Jobs and 10 Machines.

| | Tai101 | Tai102 | Tai103 | Tai104 | Tai105 | Tai106 | Tai107 | Tai108 | Tai109 | Tai110 |
|---|---|---|---|---|---|---|---|---|---|---|
| Tai101 | 0.00 | 1.00 | 1.00 | 1.00 | 1.00 | 1.00 | 1.00 | 1.00 | 0.99 | 0.99 |
| Tai102 | 1.00 | 0.00 | 1.00 | 0.99 | 1.00 | 1.00 | 0.99 | 0.97 | 0.99 | 0.99 |
| Tai103 | 1.00 | 1.00 | 0.00 | 1.00 | 0.99 | 0.98 | 1.00 | 0.99 | 1.00 | 0.99 |
| Tai104 | 1.00 | 0.99 | 1.00 | 0.00 | 1.00 | 1.00 | 1.00 | 1.00 | 1.00 | 1.00 |
| Tai105 | 1.00 | 1.00 | 0.99 | 1.00 | 0.00 | 0.98 | 1.00 | 1.00 | 0.97 | 1.00 |
| Tai106 | 1.00 | 1.00 | 0.98 | 1.00 | 0.98 | 0.00 | 1.00 | 1.00 | 0.98 | 0.98 |
| Tai107 | 1.00 | 0.99 | 1.00 | 1.00 | 1.00 | 1.00 | 0.00 | 1.00 | 0.99 | 0.98 |
| Tai108 | 1.00 | 0.97 | 0.99 | 1.00 | 1.00 | 1.00 | 1.00 | 0.00 | 1.00 | 1.00 |
| Tai109 | 0.99 | 0.99 | 1.00 | 1.00 | 0.97 | 0.98 | 0.99 | 1.00 | 0.00 | 1.00 |
| Tai110 | 0.99 | 0.99 | 0.99 | 1.00 | 1.00 | 0.98 | 0.98 | 1.00 | 1.00 | 0.00 |

Table A.13 The Inter-task Distance between Taillard's Instances with 200 Jobs and 20 Machines.

| | Tai111 | Tai112 | Tai113 | Tai114 | Tai115 | Tai116 | Tai117 | Tai118 | Tai119 | Tai120 |
|---|---|---|---|---|---|---|---|---|---|---|
| Tai111 | 0.00 | 1.00 | 1.00 | 0.99 | 0.99 | 1.00 | 1.00 | 0.98 | 1.00 | 0.99 |
| Tai112 | 1.00 | 0.00 | 0.99 | 1.00 | 1.00 | 0.99 | 0.99 | 1.00 | 0.99 | 1.00 |
| Tai113 | 1.00 | 0.99 | 0.00 | 0.99 | 1.00 | 1.00 | 0.99 | 1.00 | 1.00 | 0.99 |
| Tai114 | 0.99 | 1.00 | 0.99 | 0.00 | 1.00 | 0.99 | 1.00 | 0.99 | 1.00 | 1.00 |
| Tai115 | 0.99 | 1.00 | 1.00 | 1.00 | 0.00 | 1.00 | 0.99 | 1.00 | 1.00 | 1.00 |
| Tai116 | 1.00 | 0.99 | 1.00 | 0.99 | 1.00 | 0.00 | 0.99 | 1.00 | 1.00 | 1.00 |
| Tai117 | 1.00 | 0.99 | 0.99 | 1.00 | 0.99 | 0.99 | 0.00 | 1.00 | 0.98 | 1.00 |
| Tai118 | 0.98 | 1.00 | 1.00 | 0.99 | 1.00 | 1.00 | 1.00 | 0.00 | 0.98 | 0.98 |
| Tai119 | 1.00 | 0.99 | 1.00 | 1.00 | 1.00 | 1.00 | 0.98 | 0.98 | 0.00 | 0.99 |
| Tai120 | 0.99 | 1.00 | 0.99 | 1.00 | 1.00 | 1.00 | 1.00 | 0.98 | 0.99 | 0.00 |

Table A.14  The Inter-task Distance between Taillard's Instances with 500 Jobs and 20 Machines.